\documentclass[11pt,a4paper]{amsart}
\usepackage[dvips]{epsfig}
\usepackage{graphics}
\usepackage{latexsym}
\usepackage{verbatim}
\usepackage{amsmath}
\usepackage{amsthm}
\usepackage{amssymb}
\usepackage{MnSymbol}
\usepackage{stmaryrd}
\usepackage{hyperref}
\usepackage[]{hyperref}
\hypersetup{
    colorlinks=true,%
    filecolor=black,%
    linkcolor=black,%
    urlcolor=black  %
}
\usepackage [table]{xcolor}
\usepackage{multirow}
\usepackage{float}
\usepackage{tikz}
\usepackage{subfig}

\graphicspath{{pix/}}	
\newtheorem{theorem}{Theorem}[section]

\newtheorem{remark}[theorem]{Remark}

\usepackage{algorithmic}

\begin{document}

\title[Numerical Simulations of Kinetic Models for chemotaxis]{Numerical Simulations of Kinetic Models for chemotaxis}\thanks{The authors are  partially supported by the European Research Council ERC Starting Grant 2009,  project 239983-\textit{NuSiKiMo}}

\author{Francis Filbet and Chang Yang}

\hyphenation{bounda-ry rea-so-na-ble be-ha-vior pro-per-ties
cha-rac-te-ris-tic}

\maketitle

\begin{abstract}
We present a new algorithm based on a Cartesian mesh for the  numerical approximation of  kinetic models for chemosensitive movements set in an arbitrary geometry. We investigate the influence of the geometry on the collective behavior of bacteria described by a kinetic equation interacting with nutrients and chemoattractants. Numerical simulations are performed to verify accuracy and stability of the scheme and its ability to exhibit aggregation of cells and wave propagations. Finally some comparisons with experiments show the robustness and accuracy of such kinetic models.
\end{abstract}

\vspace{0.1cm}



\tableofcontents


\section{Introduction} 
\setcounter{equation}{0}
\label{sec:Intro}

This paper is devoted to the numerical simulation of the process of cellular spatial organization driven by chemotaxis.  The effective mechanism by which individual cells undergo directed motion varies among organisms.  Here we are particularly interested in bacterial migration, characterized by the smallness of the cells, and their ability to move to several orders of magnitude in the attractant concentration.  Several models, depending on the level of description, have been developed mathematically for the collective motion of cells. We refer to \cite{BC, ref33,ref35} for a review on parabolic, hyperbolic  and kinetic models  and to \cite{boin,bibSCBBSP, bibSCBPBS} for traveling waves drivn by growth and chemotaxis.  Among them the kinetic model introduced by  Othmer, Dunbar and Alt \cite{ref1,ref31}, describes a population of bacteria in motion (for instance the {\it E. Coli})  in  interactions with a chemoattractant concentration \cite{ref13}. These cells are so small that they are not able to measure any head-to-tail gradient of the chemical concentration, and to choose directly some preferred direction of motion towards high concentrated regions. Therefore they develop an indirect strategy to select favorable regions, by detecting a sort of time derivative in the concentration along their pathways, and to react accordingly \cite{ref27}. 

More precisely, they undergo a jump process where free runs are followed by a reorientation phenomenon called tumble \cite{ref38}. For instance it is known that {\it E. Coli} increases the time spent in running along a favourable direction \cite{ref27}. This jump process can be described by two different informations. First cells switch the rotation of their flagella, from counter-clockwise, called free runs, to clockwise called reorientation or tumbling phase, and conversely. 
This decision is the result of a complex chain of reactions inside the cells, driven by the external concentration of the chemoattractant \cite{ref37,ref38}. Then bacteria select a new direction, but they  are unable to choose directly a favourable direction, so they randomly choose a new orientation. During the "run" phases a bacterium moves with a constant speed in a given direction while during a "tumbling" event it changes direction randomly.


\section{Kinetic models for bacterial chemotaxis}
\label{sec:model}
\setcounter{equation}{0}

In the simple situation, C. S. Patlak \cite{patlak} and E. F. Keller \& L. A. Segel \cite{bibKS}  considered a density of cells which interacts with two chemical substances. The cells consume nutrients which drive the migration and  excrete a chemoattractant that prevents the dispersion of the population. However, this approach is not always sufficiently precise to describe the evolution of  bacteria movements. Hence, this phenomenon of run and tumble can be modeled by a stochastic process called the velocity-jump process, and has been introduced by Alt \cite{ref1} and further developed in \cite{ref31}. A kinetic transport model to describe this velocity jump process consists to study the evolution of the bacterial population by the local density of cells $f(t,\mathbf{x},\mathbf{v})\geq 0$ located in position $\mathbf{x}$, at time $t$ and swimming in the direction $\mathbf{v}\in V$. 

Here the set $V$ of all possible velocities is bounded and symmetric in general. In two dimensions, the modulus of the speed is a constant $v_0$, hence $V=S(0,v_0)$ circle centred in $0$ with a radius $v_0>0$. 

A kinetic transport model to describe this velocity-jump process has been introduced  by W. Alt~\cite{ref1} inspired by the Boltzmann equation~\cite{ref31} where the tumbles appear as scattering events and all the fluxes are explicitly introduced~\cite{bibSCBBSP}. We consider the Boltzmann type equation
\begin{equation}
\label{kinetic:eq}
  \frac{\partial f}{\partial t} \,+\, \mathbf{v}\cdot\nabla_\mathbf{x} f \,=\,Q(f)+r\, f, \quad  \mathbf{x}\in \Omega,\,  \quad \mathbf{v}\in V,
\end{equation} 
where $Q(f)$ is the Boltzmann type tumbling operator
\begin{equation*}
Q(f)= \int_{V} T(\mathbf{v},\mathbf{v}')\,f(t,\mathbf{x},\mathbf{v}')\,d\mathbf{v}'\,-\,\int_{V}T(\mathbf{v}',\mathbf{v})\,f(t,\mathbf{x},\mathbf{v})\,d\mathbf{v}'.
\end{equation*}
 The contribution of the tumbles is introduced with a transition (scattering) kernel $T(\mathbf{v},\mathbf{v}')\geq 0$ which stands for the change of velocity from $\mathbf{v}'$ to $\mathbf{v}$; $r$ is the division rate of the bacteria ($r=\ln2/\tau_2$ where $\tau_2$ is the mean doubling time). This equation is indeed a variant of the Boltzmann equation for gases, where collisions are delocalized via the secretion or consumption of chemical cues.

In (\ref{kinetic:eq}), the transition kernel $T$ also depends on the local concentration of chemoattractant $S(t,\mathbf{x})$ and nutrient $N(t,\mathbf{x})$. To estimate the respective contributions on pulse speed of the biais of the run lengths and of preferential reorientation, it is possible to split this transition kernel $T(\mathbf{v},\mathbf{v}')$ in two contributions, one being the tumbling rate $\lambda(\mathbf{v}')$, and the other one the reorientation effect during tumbles $K(\mathbf{v}, \mathbf{v}')$:
\begin{equation}
\label{def:T}
  T(\mathbf{v},\mathbf{v}')\,=\, \lambda(\mathbf{v}')\,K(\mathbf{v},\mathbf{v}'),
\end{equation}
with the condition 
\begin{equation}
\label{eq:Kconstant}
\int_{V}K(\mathbf{v},\mathbf{v}')\,d\mathbf{v}\,=\,1,
\end{equation}
where the function $K$ accounts for the persistence of the trajectories. For simplicity we consider the case of the absence of such an angular persistence, hence the turning kernel $T(\mathbf{v},\mathbf{v}')$ is only proportional to the tumbling rate $\lambda(\mathbf{v}')$, {\it i.e.} $K$ is constant.

For the tumbling rate $\lambda(\mathbf{v}')$, we assume that bacteria are sensitive to the temporal variations of attractant molecules via a logarithmic sensing mechanism~\cite{bibBSB, bibKJTW}. Therefore, the tumble frequency only depends on the local gradients of nutrient and attractant, both gradients having independent and additive contributions. This gives
\begin{eqnarray*}
  \lambda(\mathbf{v}')&=&\frac{1}{2}\,\left(\lambda_N(\mathbf{v}')\,+\,\lambda_S(\mathbf{v}')\right)
\\
  &=&\frac{1}{2}\,\left(\,\psi_N\left(\left.\frac{D\log N}{Dt}\right|_{\mathbf{v}'}\right)\,+\,\psi_S\left(\left.\frac{D\log S}{Dt}\right|_{\mathbf{v}'}\right)\,\right)
\\
&=&\frac{1}{2}\left(\,\psi_N\left(\partial_t \log(N)\,+\,\mathbf{v}'\cdot\nabla_\mathbf{x} \log(N)\right) \,+\,\psi_S\left(\partial_t \log(S)\,+\,\mathbf{v}'\cdot\nabla_\mathbf{x}\log(S)\right) \,\right).
\end{eqnarray*}

The nutrient and chemoattractant response functions $\psi_N$ and $\psi_S$ are both positive and decreasing, expressing that cells are less likely to tumble (thus perform longer runs) when the external chemical signal increases. These functions are smooth and characterized by their characteristic time $\delta^{-1}_N$ and $\delta^{-1}_S$ and their tumble frequency $\chi_N$ and $\chi_S$. 

Here, we have chosen the following analytical form that encompassed these characteristics:
\begin{equation}
\label{def:psi}
  \psi_\alpha(X)=\psi_0\left(1-\chi_\alpha\tanh\left(\frac{X}{\delta_\alpha}\right)\right), \quad \alpha=\{N,S\},
\end{equation}
where $\psi_0$ is the mean tumbling frequency, and the parameters $\chi_N$ and $\chi_S$ are the modulation of tumble frequency.

In order to define completely the mathematical problem \eqref{kinetic:eq}, suitable boundary conditions on $\partial\Omega$ should be applied. Here we consider wall type boundary conditions, for which  emerging  particles have been reflected elastically at the wall. More precisely, for $\mathbf{x}\in\partial\Omega$ the smooth boundary $\partial\Omega$ is assumed to have a unit inward normal $\mathbf{n}(\mathbf{x})$ and for $\mathbf{v}\cdot\mathbf{n}(\mathbf{x})\geq0$, we assume that at the solid boundary we have
\begin{equation}
  f(t,\mathbf{x},\mathbf{v})\,=\,\mathcal{R}[f(t,\mathbf{x},\mathbf{v})],\,\,\,\mathbf{x}\in\partial\Omega,\,\,\,\mathbf{v}\cdot\mathbf{n}(\mathbf{x})\geq 0,
\label{eq:BCingoing}
\end{equation}
with
\begin{equation}
\label{op:BC}
\mathcal{R}[f(t,\mathbf{x},\mathbf{v})] = f(t,\mathbf{x},\mathbf{v}-2(\mathbf{v}\cdot\mathbf{n}(\mathbf{x}))\mathbf{n}(\mathbf{x})).
\end{equation}
This boundary condition~\eqref{eq:BCingoing} guarantees the global conservation of mass~\cite{bibC}.

The equations describing the behaviors of nutrients density $N$ and chemottractant $S$ are the same as in~\cite{bibMBBO}:
\begin{equation}
\label{NS:eq}
  \left\{
  \begin{array}{l}
    \displaystyle \frac{\partial S}{\partial t} \,-\,D_S\,\Delta S \,=\, -\,a \,S\,+\,b\,\int_{V} f(t,\mathbf{x},\mathbf{v})\,d\mathbf{v}, \quad \mathbf{x}\in \Omega,
\\
\,
\\
     \displaystyle \frac{\partial N}{\partial t} \,-\,D_N\,\Delta N\,=\,-\,c\,N\,\int_{V} f(t,\mathbf{x},\mathbf{v})\,d\mathbf{v}, \quad \mathbf{x}\in \Omega,
  \end{array}
  \right.
\end{equation}
where $a$, $b$ and $c$ are respectively the degradation rate of the chemoattractant, its production rate and the consumption rate of the nutrient by the bacteria, whereas $D_S$ and $D_N$ are the molecular diffusion coefficients. Finally, these equations are completed with homogeneous Neumann boundary conditions, {\it i.e.} 
\begin{equation}
  \nabla_{\mathbf{x}}\alpha\cdot\mathbf{\mathbf{n}(\mathbf{x})}=0,\quad \alpha=\{N,S\},\quad\mathbf{x}\in \partial\Omega.
\label{eq:BC:parabolique}
\end{equation}

The purpose of this work is to present  a numerical scheme for  (\ref{kinetic:eq})-(\ref{eq:BC:parabolique}) and to investigate numerically the occurrence of cells aggregation, pattern formation or travelling waves when it takes place, and the convergence to equilibrium otherwise for different geometries.  Several numerical methods have already been developed to solve the Patlak-Keller-Segel model for chemotaxis using finite element methods \cite{marrocco},  finite volume methods \cite{Chertock, ref:ep, filbet:PKS}, and the references therein. Other models have also been investigated numerically \cite{filbet1,filbet2,ribot}. However, it seems that none of the above-mentioned numerical approaches have been studied for kinetic models (\ref{kinetic:eq})-(\ref{eq:BC:parabolique}). In the present paper we propose a numerical scheme for (\ref{kinetic:eq})-(\ref{eq:BC:parabolique})  and investigate the influence of the geometry on the collective behavior of bacteria.   

We now briefly outline the contents of the paper. In the next section, we introduce the numerical approximation of (\ref{kinetic:eq}) and (\ref{NS:eq}) and describe the numerical approximation of the boundary~\eqref{eq:BCingoing},~\eqref{op:BC},~\eqref{eq:BC:parabolique}. Two points are worth mentioning here. First, we restrict ourselves to the case of specular reflection which seems to be the most appropriate for the study of bacteria. One difficulty in the approximation of kinetic models for chemotaxis, is related to the fact that it can exhibits very different phenomena as finite time blow-up, cell aggregation, wave propagations. At the discrete level, our approximation should also be able to describe  a similar property. Secondly, an important step is to discretize appropriately the effect of  boundary.

 The final section is devoted to numerical simulations performed with the numerical scheme presented in Section \ref{sec:scheme}. We investigate numerically cells aggregation, convergence to equilibrium, and wave propagation  in a bounded domain.


\section{Numerical resolution}
\label{sec:scheme}
\setcounter{equation}{0}


\subsection{Numerical resolution of the kinetic model (\ref{kinetic:eq})}
We consider a computational domain $[x_{\min},x_{\max}]\times[y_{\min},y_{\max}]\times V$,  such that $\Omega\subset[x_{\min},x_{\max}]\times[y_{\min},y_{\max}]$. 

The computational domain is covered by an uniform Cartesian mesh $\mathbf{X}_h\times \mathbf{V}_{h}$, where $\mathbf{X}_h$, $\mathbf{ V}_{h}$ are defined by 
\begin{equation}
\left\{
\begin{array}{l}
  \mathbf{X}_h\,:=\,\left\{\,\mathbf{x}_{0}=(x_{\min},y_{\min}),\ldots,  \mathbf{x}_{i}=(x_{i_x},y_{i_y}),\ldots, \mathbf{x}_{(n_x,n_y)}=(x_{\max},y_{\max})\,   \right\},
\\
\,
\\
 \mathbf{V}_h\,:=\,\left\{\mathbf{v}_j\,=\, v_0\,(\cos\theta_j,\sin\theta_j),\quad  \theta_j = (j+1/2)\,\Delta v, 0\leq j\leq n_v-1 \right\}.
\end{array}
\right.
\label{eq:2Dmesh}
\end{equation}
The mesh steps are respectively $\Delta x=(x_{\max}-x_{\min})/n_x$, $\Delta y=(y_{\max}-y_{\min})/n_y$ and $\Delta v=2\pi/n_v$.

Let us denote $f^n_{i,j}$ an approximation of the distribution function $f(t^n,\mathbf{x}_{i},\mathbf{v}_j)$. We introduce  the following finite difference scheme
\begin{equation}
\label{kinetic:eq:time}
  \frac{f^{n+1}_{i,j}-f^n_{i,j}}{\Delta t} \,+\, \mathbf{v}\cdot\nabla_{\mathbf{x},h} f^n_{i,j} \,=\, Q_h(f^n_{i,j})+r\cdot f^n_{i,j},
\end{equation} 
where $h=(\Delta t, \Delta x,\Delta y)$, $\mathbf{v}\cdot\nabla_{\mathbf{x},h} f^n_{i,j}$ is a second-order approximation~\cite{bibVL} of the transport operator $\mathbf{v}\cdot\nabla_\mathbf{x} f$, and $Q_h(f^n_{i,j})$ is an approximation of the Boltzmann type tumbling operator $Q(f)$. We will now focus on searching the approximation $Q_h(f^n_{i,j})$.

By using the trapezoidal rule, we have
\begin{equation*}
  Q_h(f^n_{i,j})=\Delta v\sum_{\ell=0}^{n_v-1}\left((T_h)^{n+1/2}_{i,j,\ell}f^n_{i,\ell}\right)-\Delta v\sum_{\ell=0}^{n_v-1}\left((T_h)^{n+1/2}_{i,\ell,j}f^n_{i,j}\right),
\end{equation*}
where $(T_h)^{n+1/2}_{i,j,\ell}$ is an approximation of the transition kernel $T(\mathbf{v},\mathbf{v}')$. 
 We assume that the nutrient density $N$ and the chemottractant $S$ are known at time $t_{n+1/2}$. Moreover with the hypothesis that $K$ is constant, the condition~\eqref{eq:Kconstant} implies that $K={1}/{2\pi}$. Thus it remains to search an approximation of the tumbling rate, {\it i.e.} $\lambda^{n+1/2}_{i,\ell}$. It is also equivalent to search the local gradients of nutrient  $(\lambda_N)^{n+1/2}_{i,\ell}$ and attractant $(\lambda_S)^{n+1/2}_{i,\ell}$.  

We  study only the local gradients of nutrient $(\lambda_N)^{n+1/2}_{i,\ell}$, since the local gradient of  attractant $(\lambda_S)^{n+1/2}_{i,\ell}$ has the same expression as $(\lambda_N)^{n+1/2}_{i,\ell}$. We discretize the local gradient of nutrient as follows
\begin{equation*}
  (\lambda_N)^{n+1/2}_{i,\ell}=\frac{1}{N^{n+1/2}_{i}}\left(D_{h} N^{n+1/2}_{i}+\mathbf{v}_{\ell}\cdot\nabla_{\mathbf{x},h}N^{n+1/2}_{i}\right),
\end{equation*}
where  $D_{h} N^{n+1/2}_{i}$ is a discrete time derivative and will be given in the  section~\ref{sec:discratisation:parabolique}. Moreover we use  centred difference approximation for $\mathbf{v}'\cdot\nabla_{\mathbf{x}}N$, which yields
\begin{equation*}
  \mathbf{v}_{\ell}\cdot\nabla_{\mathbf{x},h}N^{n+1/2}_{i}=v_0\cos\theta_{\ell}\frac{N^{n+1/2}_{i_x+1,i_y}-N^{n+1/2}_{i_x-1,i_y}}{2\Delta x}+v_0\sin\theta_{\ell}\frac{N^{n+1/2}_{i_x,i_y+1}-N^{n+1/2}_{i_x,i_y-1}}{2\Delta y}.
\end{equation*}
 In summary, the discrete tumbling rate reads

\begin{equation*}
  \lambda^{n+1/2}_{i,\ell}=\frac{1}{2}\left((\lambda_N)^{n+1/2}_{i,\ell}+(\lambda_S)^{n+1/2}_{i,\ell}\right).
\end{equation*}
Finally we reduce the tumbling operator as follows
\begin{equation*}
  Q_h(f^n_{i,j})=\frac{\Delta v}{2\pi}\sum_{\ell=0}^{n_v-1}\left(\lambda^{n+1/2}_{i,\ell}f^n_{i,\ell}\right)-\lambda^{n+1/2}_{i,j}f^n_{i,j}.
\end{equation*}

\subsection{Calculate the chemoattractant $S$ or the nutrient $N$}
\label{sec:discratisation:parabolique}

The  equations~\eqref{NS:eq} for nutrient density $N$ and chemottractant $S$ are parabolic equations with source terms depending on the distribution function of density $f$. We study again the discretization for nutrient density, since the one for chemottractant is similar. The Euler implicit scheme is used for time integration. Hence the scheme for nutrient density $N$ reads
\begin{equation}
  D_{h} N^{n+1/2}_{i}-D_N\Delta_hN_i^{n+1/2}=-cN^{n+1/2}_i\rho_i^n,
\end{equation}
where $ D_h N^{n+1/2}_{i}$ is an approximation of time derivative as follows
\begin{equation}
  D_h N^{n+1/2}_{i}=\frac{N^{n+1/2}_{i}-N^{n-1/2}_{i}}{\Delta t}.
\end{equation}
Then we use a five points finite difference scheme to discretize $\Delta N$
\begin{equation}
\label{eq:discrete:laplace}
  \Delta_h N^{n+1/2}_{i}=\frac{N^{n+1/2}_{i_x+1,i_y}-2N^{n+1/2}_{i_x,i_y}+N^{n+1/2}_{i_x-1,i_y}}{\Delta x^2}+\frac{N^{n+1/2}_{i_x,i_y+1}-2N^{n+1/2}_{i_x,i_y}+N^{n+1/2}_{i_x,i_y-1}}{\Delta y^2}.
\end{equation}
Finally, a trapezoidal rule is applied for the integration, {\it i.e.}
\begin{equation*}
 \rho^n_i=\Delta v\sum_{j=0}^{n_v-1}f^n_{i,j}.
\end{equation*}

\subsection{Treatment of the boundary conditions}
\label{sec:discret_BC}
As we mentioned at the end of section~\ref{sec:model}, an appropriate discretization of boundary condition is important to exhibit very different phenomena. Therefore, we  present respectively the numerical approximations for specular reflection boundary condition~\eqref{op:BC} and Neumann boundary condition~\eqref{eq:BC:parabolique}.
\subsubsection{Numerical approximation for specular reflection boundary condition}


 The specular reflection boundary condition in 2D reads as
 \begin{equation}
   \mathcal{R}[f](t,\mathbf{x},\mathbf{v})=f(t,\mathbf{x},\mathbf{v}'),
 \end{equation}
 with
 \begin{equation*}
  \mathbf{v}\cdot\mathbf{n}=-\mathbf{v}'\cdot\mathbf{n} \,\,\Leftrightarrow\,\,\mathbf{v}'=\mathbf{v}-2( \mathbf{v}\cdot\mathbf{n})\mathbf{n}
 \end{equation*}
where $\mathbf{x}\in\Gamma_{\mathbf{x}}=\partial\Omega$ is the point at the boundary, $\mathbf{n}$ is the interior normal at point $\mathbf{x}$. We note that this specular reflection is just like a mirror. For example, if we follows the characteristic of the flux $f(\mathbf{v})$, its reflected flow is corresponding to the velocity $\mathbf{v}'$.

\begin{figure}[h]
  \begin{center} 
    \includegraphics[width=10cm]{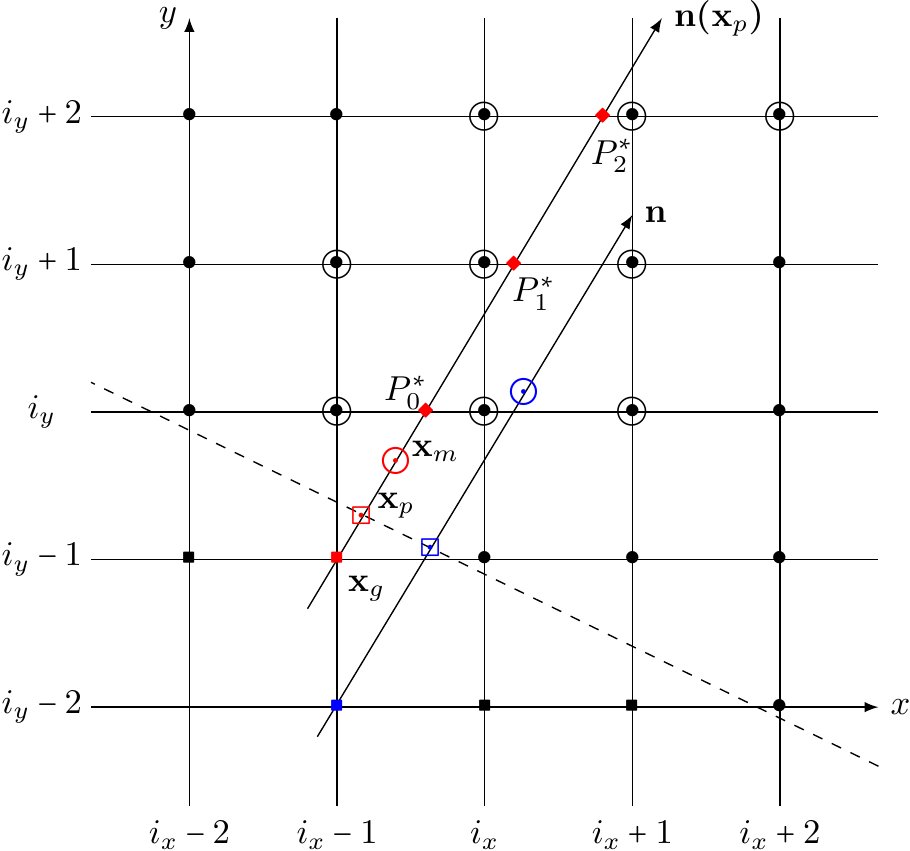}
\caption{\label{fig:2Ddomain}Spatially two-dimensional Cartesian mesh. $\bullet$ is interior point, $\filledsquare$ is ghost point, $\boxdot$ is the point at the boundary, $\largecircle$ is the point for extrapolation, the dashed line is the boundary.}
  \end{center}
\end{figure}

We thus use a mirror procedure to construct $f$ at each ghost point. For instance from the ghost point $\mathbf{x}_g$, we can find an inward normal $\mathbf{n}(\mathbf{x}_p)$, which  crosses the boundary at $\mathbf{x}_p$ (see Figure~\ref{fig:2Ddomain}). For velocity $\mathbf{v}$, its reflected velocity with respect to $\mathbf{x}_p$ is $\mathbf{v}'=\mathbf{v}-2( \mathbf{v}\cdot\mathbf{n}(\mathbf{x}_p))\mathbf{n}(\mathbf{x}_p)$. Thus instead of computing  $f$ at the ghost point $\mathbf{x}_g$, we compute $f$ at  mirror point with respect to the boundary $\mathbf{x}_m=2\mathbf{x}_p-\mathbf{x}_g$ as follows
\begin{equation}
  f(\mathbf{x}_g,\mathbf{v})=f(\mathbf{x}_m,\mathbf{v}').
\end{equation}

The last step is to approximate $f(\mathbf{x}_m,\mathbf{v}')$ using $f$ of interior domain. Let us assume that the values of the distribution function $f$ on the grid points in $\Omega$ are given. We first construct a stencil $\mathcal{E}$ composed of grid points of $\Omega$ for interpolation or extrapolation. For instance as it is shown in Figure~\ref{fig:2Ddomain},  the inward normal $\mathbf{n}(\mathbf{x}_p)$ intersects the grid lines $y=y_{i_y}$, $y_{i_y+1}$, $y_{i_y+2}$ at points $P^*_0$, $P^*_1$, $P^*_2$. Then we choose the three nearest points  of the cross point $P^*_l,\,l=0,1,2$, in each line, {\it i.e.} marked by a large circle. From these nine points, we can build a  Lagrange polynomial $q_2(\mathbf{x})\in\mathbb{Q}_2(\mathbb{R}^2)$. Therefore  we evaluate the polynomial $q_2(\mathbf{x})$ at $\mathbf{x}_m$, and obtain an approximation of $f$ at the mirror point. We distinguish two cases of mirror points. In the case that mirror point $\mathbf{x}_m$ is surrounded by interior points, we interpolate $f$ at mirror point  $\mathbf{x}_m$; otherwise a WENO type extrapolation can be used to prevent spurious oscillations, which will be detailed below.

Note that in some cases, we can not find a stencil of nine interior points. For instance when the interior domain has small acute angle sharp, the normal $\mathbf{n}$ can not have three cross points $P^*_l,\,l=0,1,2$ in interior domain, or we can not have three nearest points  of the cross point $P^*_l,\,l=0,1,2$, in each line. In such a case, we alternatively use a first degree polynomial $q_1(\mathbf{x})$ with a four points  stencil or even a zero degree polynomial $q_0(\mathbf{x})$ with  an one point  stencil. We can similarly construct the four points stencil or the one point stencil as above.

\paragraph{\bf Two-dimensional WENO type extrapolation}A WENO type extrapolation~\cite{bibTS} was developed to prevent oscillations and maintain accuracy, which is an extension of WENO scheme~\cite{bibJS}. The key point of WENO type extrapolation is to define smoothness indicators, which is designed to help us  choose automatically between the high order accuracy  and the low order but more robust extrapolation. Moreover  a slightly modified version of the  method was given in~\cite{bibFY}, such that  the smoothness indicators are invariant with respect to the scaling of $f$. We  now describe this method in 2D  case. 

The substencils $S_r,\,\,r={0,1,2}$ for extrapolation are chosen around the inward normal $\mathbf{n}$ such that we can construct Lagrange polynomial of degree $r$. For instance in  Figure~\ref{fig:2Ddomain}, the three substencils are respectively  
$$
\left\{\begin{array}{l}
 \displaystyle S_0\,=\;\{(x_{i_x},y_{i_y})\},
\\ \,\\
 \displaystyle S_1\,=\,\{(x_{i_x-1},y_{i_y}),(x_{i_x},y_{i_y}),(x_{i_x},y_{i_y+1}),(x_{i_x+1},y_{i_y+1})\},
\\\,\\
 \displaystyle S_2\,=\,  \{(x_{i_x-1},y_{i_y}),(x_{i_x},y_{i_y}),(x_{i_x+1},y_{i_y}),(x_{i_x-1},y_{i_y+1}),
\\\,\\
 \displaystyle\,\quad (x_{i_x},y_{i_y+1}),(x_{i_x+1},y_{i_y+1}),(x_{i_x},y_{i_y+2}),(x_{i_x+1},y_{i_y+2}),(x_{i_x+2},y_{i_y+2})\}.
\end{array}\right.
$$
Once the substencils $S_r$ are chosen, we could easily construct the Lagrange polynomials in $\mathbb{Q}_r(\mathbb{R}^2)$
\begin{equation*}
  q_r(\mathbf{x})=\sum_{m=0}^r\sum_{l=0}^ra_{l,m}x^ly^m
\end{equation*}
satisfying 
\begin{equation*}
  q_r(\mathbf{x})=f(\mathbf{x}),\,\,\mathbf{x}\in S_r.
\end{equation*}
Then the  WENO extrapolation has the form
\begin{equation}
  f(\mathbf{x})=\sum_{r=0}^2w_rq_r(\mathbf{x}),\,\,\mathbf{x}\in S_r,
\end{equation}
where $w_r$ are the nonlinear weights, which are chosen to be
\begin{equation*}
  w_r=\frac{\alpha_r}{\sum_{s=0}^2\alpha_s},
\end{equation*}
with 
\begin{equation*}
  \alpha_r=\frac{d_r}{(\varepsilon+\beta_r)^2},
\end{equation*}
where $\varepsilon=10^{-6}$, $d_0=\Delta x^2+\Delta y^2$, $d_1=\sqrt{\Delta x^2+\Delta y^2}$, $d_2=1-d_0-d_1$. The coefficients $\beta_r$ are the smoothness indicators determined by
\begin{eqnarray*}
  \beta_0&=&\Delta x^2+\Delta y^2,
\end{eqnarray*}
and for $r\geq 1$, either we take $\beta_r=0$  when $f(\mathbf{x})=0,\forall\mathbf{x}\in S_r$, or  we choose
\begin{eqnarray*}
  \beta_r&=&\frac{1}{\sum_{\mathbf{x}\in S_r}f(\mathbf{x})^2}\sum_{1\leq|\sigma|\leq r}\int_K|K|^{|\sigma|-1}(D^\sigma q_r(\mathbf{x}))^2d\mathbf{x},\,\,r=1,2,
\end{eqnarray*}
where $\sigma$ is a multi-index and $K=[x_p-\Delta x/2,x_p+\Delta x/2]\times[y_p-\Delta y/2,y_p+\Delta y/2]$ and the point  $\mathbf{x}_p$ is given by $(x_p,y_p)$.

\subsubsection{Numerical approximation for Neumann boundary condition}
\begin{figure}[h]
  \begin{center} 
    \includegraphics[width=10cm]{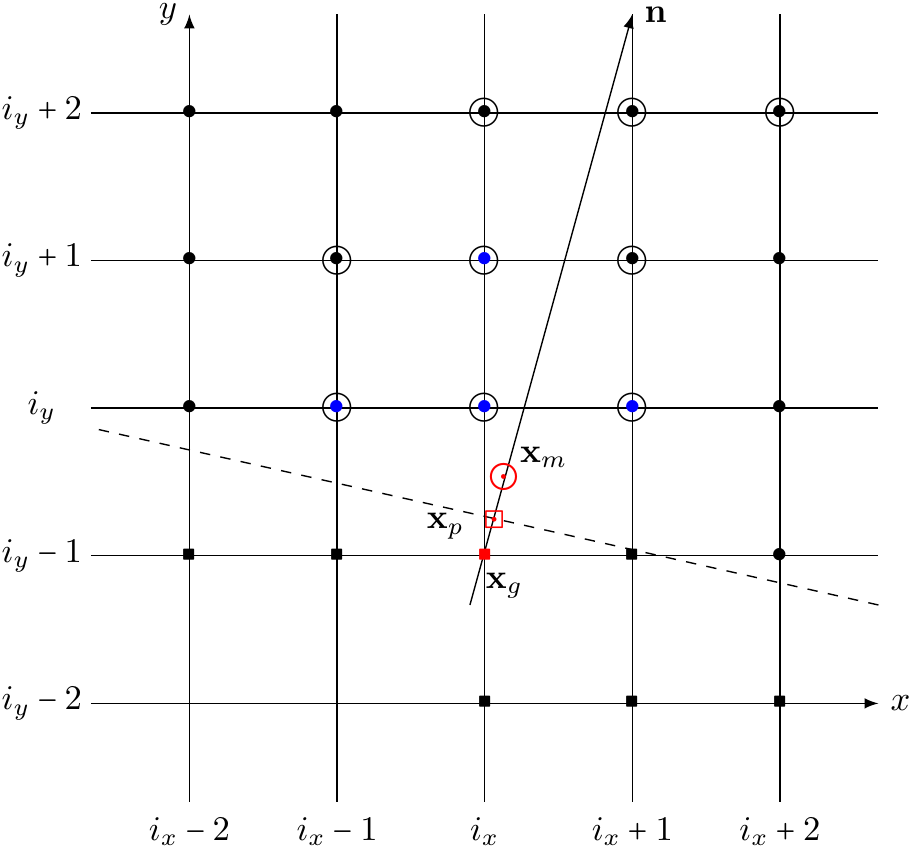}
\caption{\label{fig:2Ddomain_parapolique}Illustration of discretization of Laplace operator near the boundary on two-dimensional Cartesian mesh. $\bullet$ is interior point, $\filledsquare$ is ghost point, $\boxdot$ is the point at the boundary, $\largecircle$ is the point for extrapolation, the dashed line is the boundary.}
  \end{center}
\end{figure}

In section~\ref{sec:discratisation:parabolique}, we have seen that discrete Laplace operator~\eqref{eq:discrete:laplace} is the only non-locally in space term in~\eqref{NS:eq}. In Figure~\ref{fig:2Ddomain_parapolique}, we illustrate an example of a discrete  Laplace operator near the boundary. We note that $N_{i_x,i_y-1}$ is not known, since $\mathbf{x}_g=(x_{i_x},y_{i_y-1})$ is out of domain. To approximate $N$ at ghost point $\mathbf{x}_g$, we  have to use the boundary condition~\eqref{eq:BC:parabolique}.

In fact, if we denote $\mathbf{n}:=(n_x,n_y)$, then the boundary condition~\eqref{eq:BC:parabolique} reads 
\begin{equation*}
  n_x\partial_x N+n_y\partial_y N=0.
\end{equation*}
Using a centered difference formula, it yields
\begin{equation*}
  n_x\frac{N(\mathbf{x}_m)-N(\mathbf{x}_g)}{x_m-x_g}+n_y\frac{N(\mathbf{x}_m)-N(\mathbf{x}_g)}{y_m-y_g}=0,
\end{equation*}
where $\mathbf{x}_g:=(x_g,y_g)$, $\mathbf{x}_m:=(x_m,y_m)$ is the mirror point of the ghost point $\mathbf{x}_g$ with respect to the boundary. The previous equation is equivalent to 
\begin{equation}
  N(\mathbf{x}_m)=N(\mathbf{x}_g).
\end{equation}
Therefore instead of computing $N$ at the ghost point $\mathbf{x}_g$, we extrapolate $N$ at the mirror point  $\mathbf{x}_m:=(x_m,y_m)$. As shown in Figure~\ref{fig:2Ddomain_parapolique}, a nine points stencil $S_2$ is found, thus we have
\begin{equation}
  \label{eq:mirror:discret}
  N_{i_x,i_y-1}=N(\mathbf{x}_m)=\sum_{i=0}^{8}w_iN(\mathbf{x}_i),\quad\,\mathbf{x}_i\in\,S_2,
\end{equation}
where $w_i$ is weight of extrapolation calculated, for instance, by Lagrange polynomial. Finally by injecting~\eqref{eq:mirror:discret} into~\eqref{eq:discrete:laplace}, we complete the discretization.

\begin{remark}
It is important to emphasize that the numerical scheme for the internal domain is disconnected to the numerical procedure for the treatment of the boundary. Therefore the present method can be applied to any other numerical scheme. The main point is to preserve the order of accuracy and then to eventually increase the stencil outside the domain.
\end{remark}

\section{Numerical examples}
\label{sec:Num}
\setcounter{equation}{0}
In this section, we present a large variety of tests in $2d_x$ and one dimensional in velocity space showing the validation of mathematical model and the effectiveness of our numerical schemes. In Test 1, we consider cell aggregation to study only chemotactic motility~\cite{bibMBBO}. In Test 2, we consider  wave propagation  formed by cells reorientation and the presence of nutrients in a disc. In Test 3, we study the interaction of two traveling wave in a $U$ shape. Then in Test 4, we consider again the traveling wave but in a more wide $U$ shape, which is similar like the effectiveness of centrifugal force. Finally in Test 5 we add the effect of cell division for long time behavior using new parameters.


\subsection{Test 1 : cell aggregation}
Motile bacteria are able to interact with their environment by accumulating in regions of high concentrations of certain chemicals called attractants and avoiding others called repellents. Cell division is not required to generate these multi-cellular structures because they form on a much shorter time scale than the cell-doubling time~\cite{bibMBBO}. The process leading to formation of multi-cellular clusters can be qualitatively understood as follows: fluctuations in the local cell density produce local gradients of attractants. Cells respond by moving up these concentration gradients thus amplifying the initial spatial non-uniformities in the cell distribution and  forming multi-cellular clusters.

We thus use the model of Othmer-Dunbar-Alt~\cite{ref31} to mimic the cell aggregation as follows
\begin{equation*}
  \partial_t f+\mathbf{v}\cdot\nabla_{\mathbf{x}}f=\int_{\mathbf{v}'}T[S](\mathbf{v},\mathbf{v}')f(t,\mathbf{x},\mathbf{v}')d\mathbf{v}'-\int_{\mathbf{v}'}T[S](\mathbf{v}',\mathbf{v})f(t,\mathbf{x},\mathbf{v})d\mathbf{v}',
\end{equation*}
where the tumbling kernel $T[S](\mathbf{v}',\mathbf{v})$ describes the frequency of reorientation $\mathbf{v}'\to\mathbf{v}$ 
\begin{equation*}
T[S](\mathbf{v}',\mathbf{v})=\psi_S\left(\partial_t \log(S)\,+\,\mathbf{v}'\cdot\nabla_\mathbf{x}\log(S)\right) ,
\end{equation*}
and the chemical signal is secreted by the cells, following the reaction-diffusion equation~\eqref{NS:eq}.

We use a square domain of size $[-0.25,0.25]^2$, which is uniformly divided by a mesh size of $n_x\times n_y=120\times120$. The velocity space belongs to  the unit circle $S^1$, and is uniformly divided into $n_v=64$ parts. All the parameters are chosen as in~\cite{bibSCBPBS} and listed in Table~\ref{tab:parameters}. The initial chemoattractant $S_0$ is equal to 0, and the initial distribution function $f_0$ is independent of the velocity $\mathbf{v}$
\begin{equation*}
  f_0(\mathbf{x},\mathbf{v})=\frac{100m}{\pi}\exp(-100|\mathbf{x}|^2),
\end{equation*}
where $m$ is the total mass.

\begin{table}
  \caption{\label{tab:parameters}Summary of the values used in the simulation}
  \begin{tabular}{|l|l|}
\hline
    Parameter&Value\\
    \hline
    run speed&$v_0=25\mu m\cdot s^{-1}$\\
    mean tumbling frequency&$\psi_0=3 s^{-1}$\\
    modulation of tumbling frequency of nutrient&$\chi_N=60\%$\\
    modulation of tumbling frequency of chemoattractant&$\chi_S=20\%$\\
    stiffness of the response functions&$1/\delta_N=1/\delta_N\approx20s$\\
    space scale&$\bar x=1mm$\\
    time scale&$\bar t=40s$\\
    doubling time&$\tau_2=\ln2/r\approx2h$\\
    degradation rate of the chemoattractant&$a=5\cdot10^{-3}mol\cdot s^{-1}$\\
    production rate of the chemoattractant&$b=4\cdot10^5cell^{-1}s^{-1}$\\
    consumption rate of the nutrient by the bacteria&$c=2\cdot10^{-7}cell^{-1}s^{-1}$\\
    diffusion coefficient of the nutrient molecules&$D_N=8\cdot10^{-6}cm^2\cdot s{-1}$\\
    diffusion coefficient of the chemoattractant molecules&$D_S=8\cdot10^{-6}cm^2\cdot s{-1}$\\
    
    \hline
  \end{tabular}
\end{table}

\begin{figure}[htbp]
\begin{tabular}{ccc}
\includegraphics[width=4.2cm]{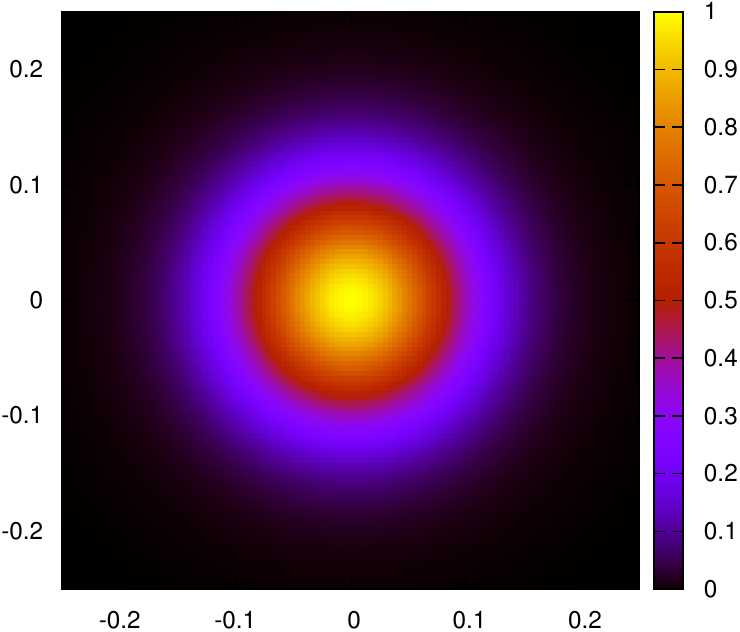} &
\includegraphics[width=4.5cm]{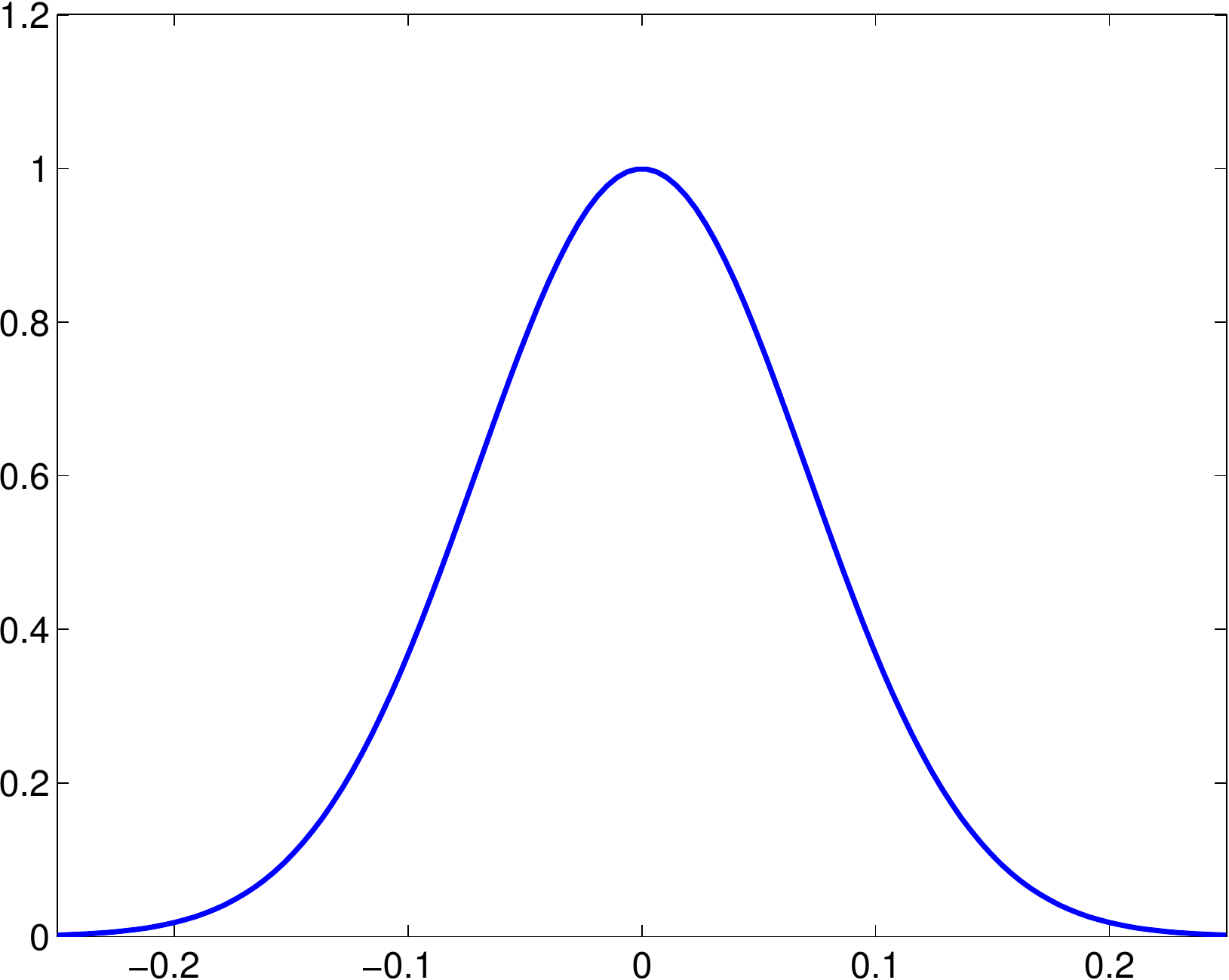} &
\includegraphics[width=4.5cm]{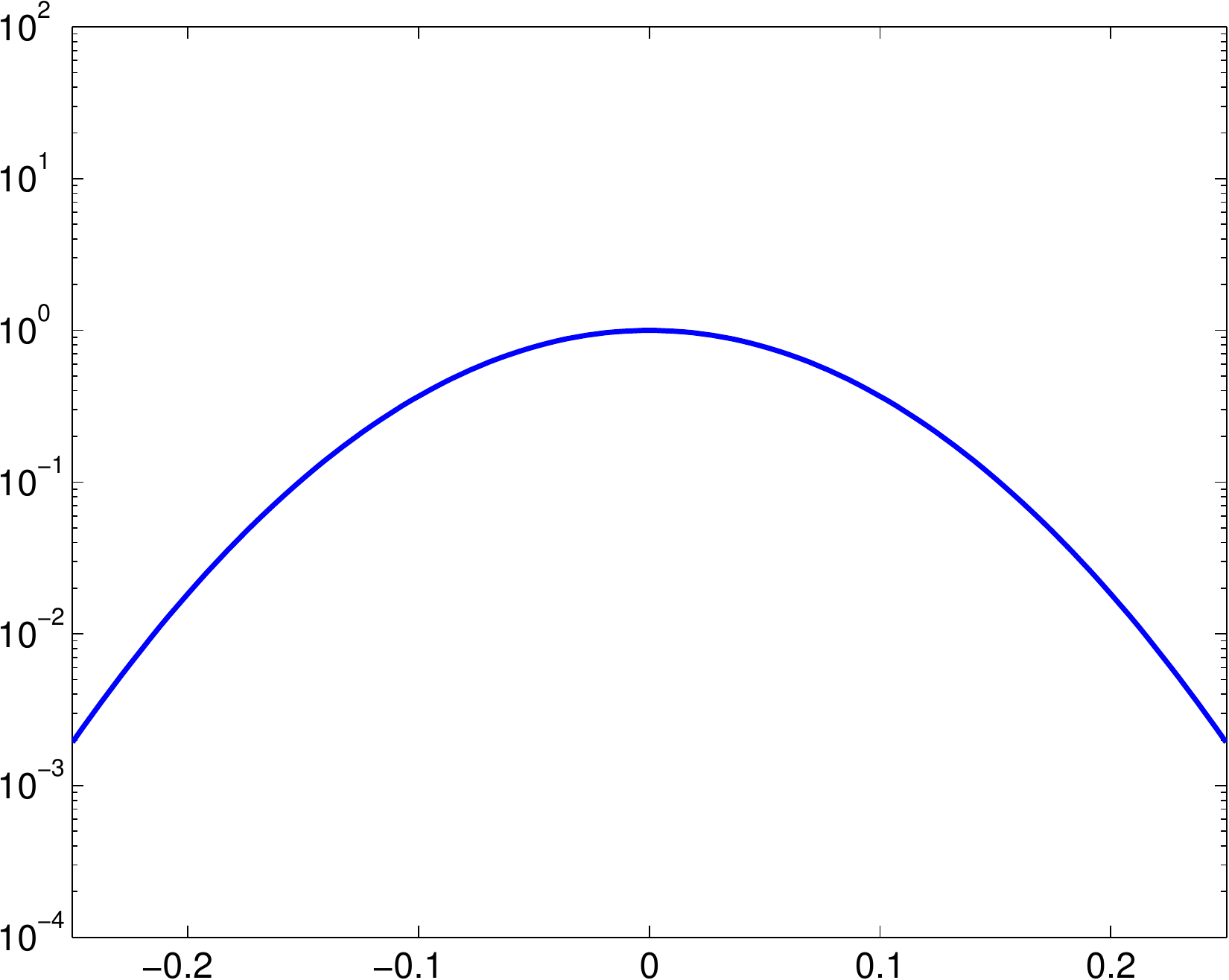} 
\\
$t=0$& $t=0$& $t=0$
\\
\includegraphics[width=4.2cm]{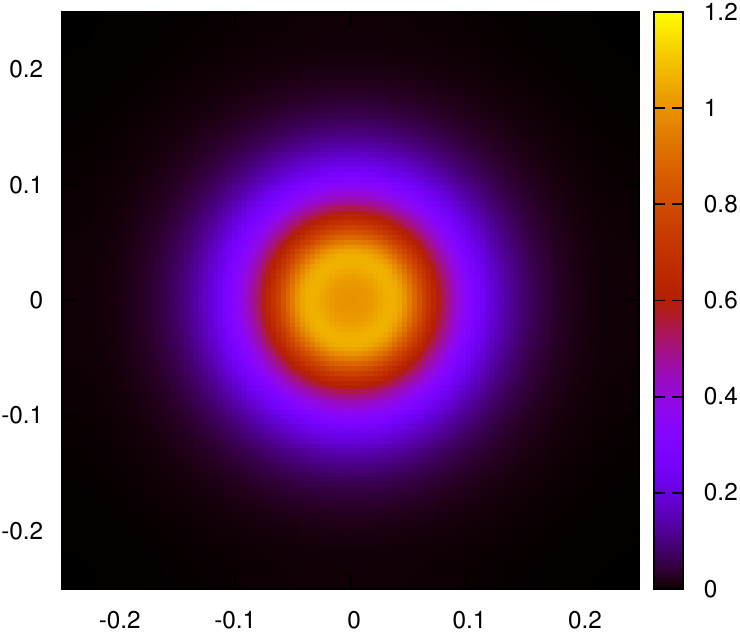} &
\includegraphics[width=4.5cm]{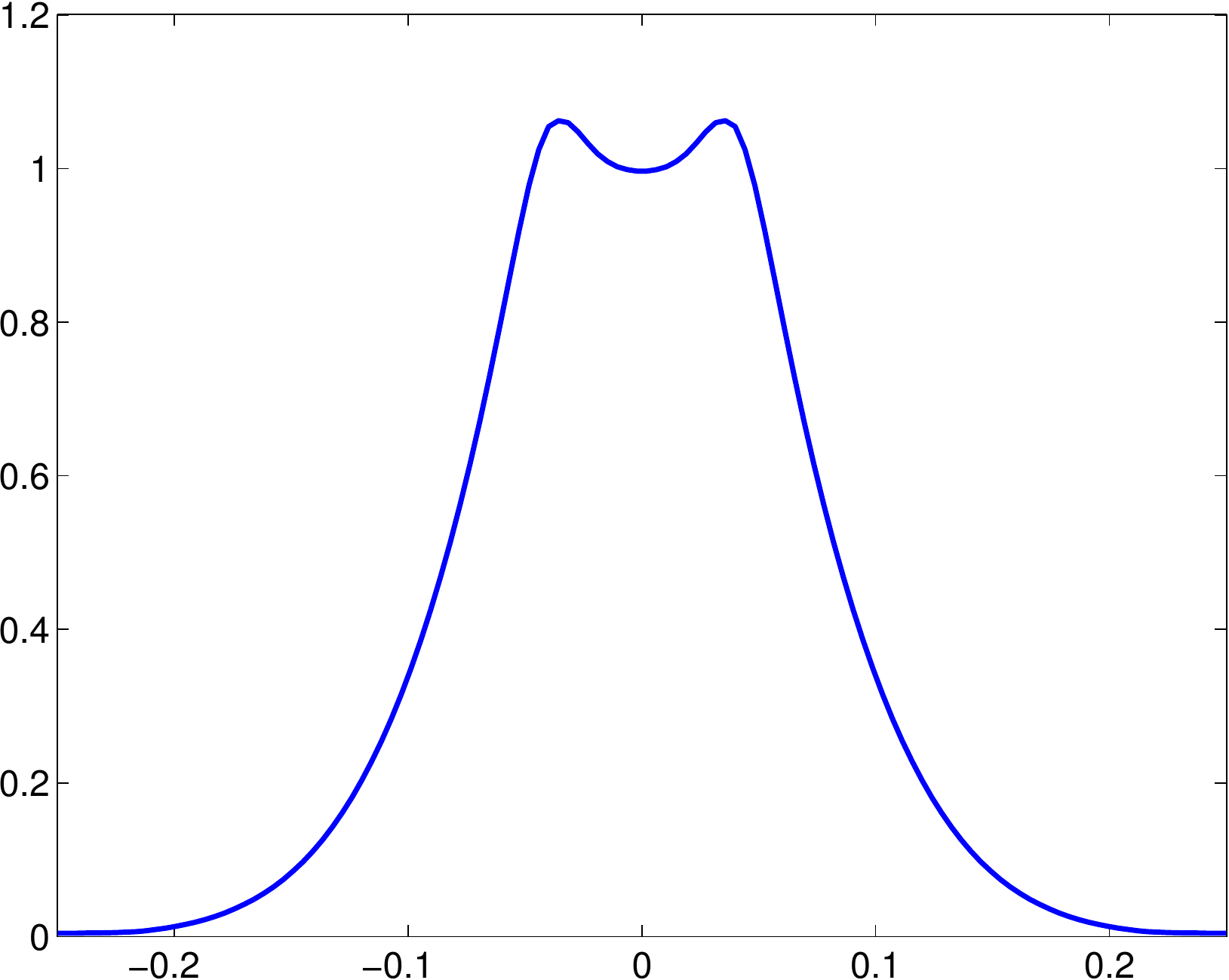} &
\includegraphics[width=4.5cm]{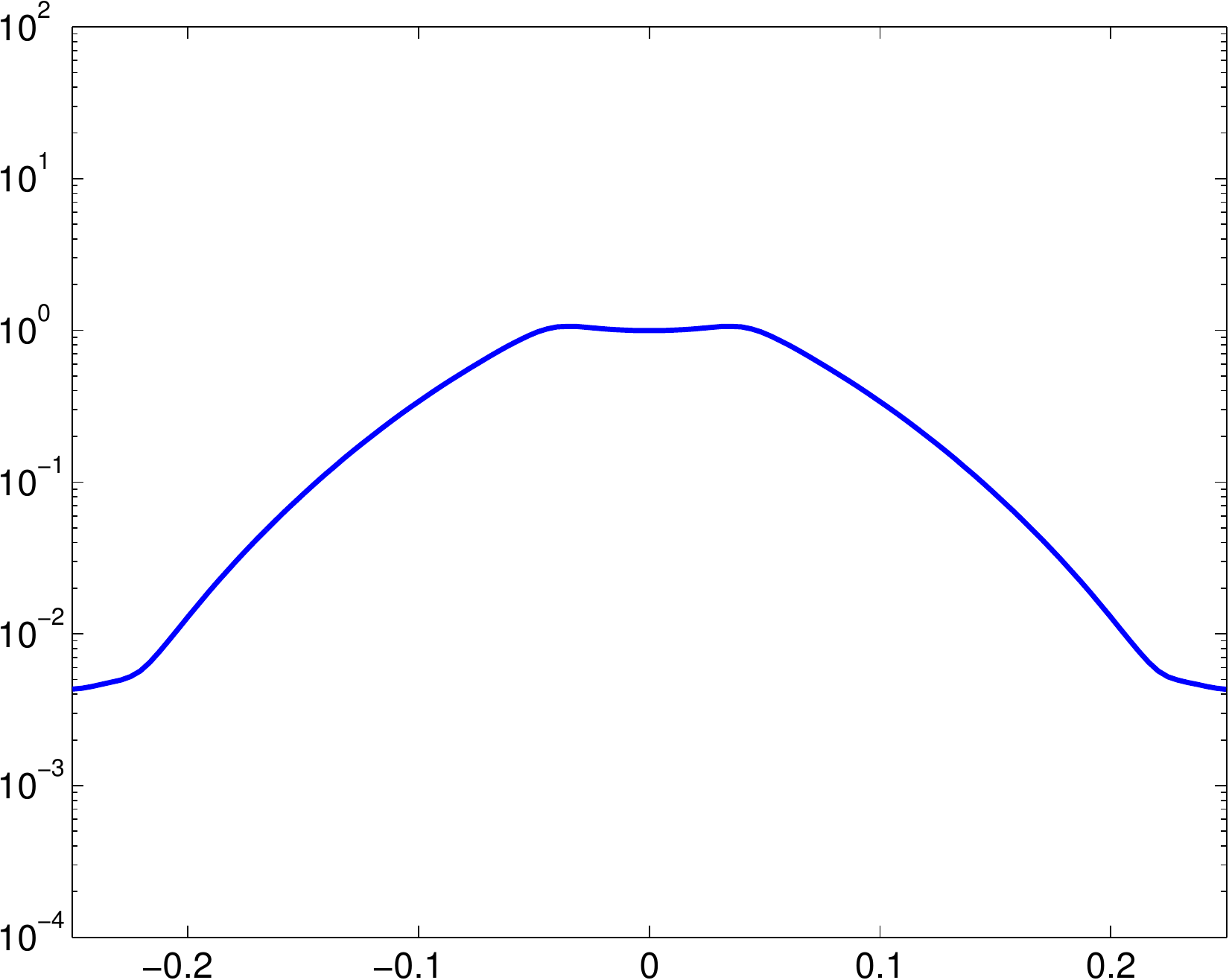} 
\\
$t=0.2\,\bar t$& $t=0.2\,\bar t$& $t=0.2\,\bar t$
\\
\includegraphics[width=4.2cm]{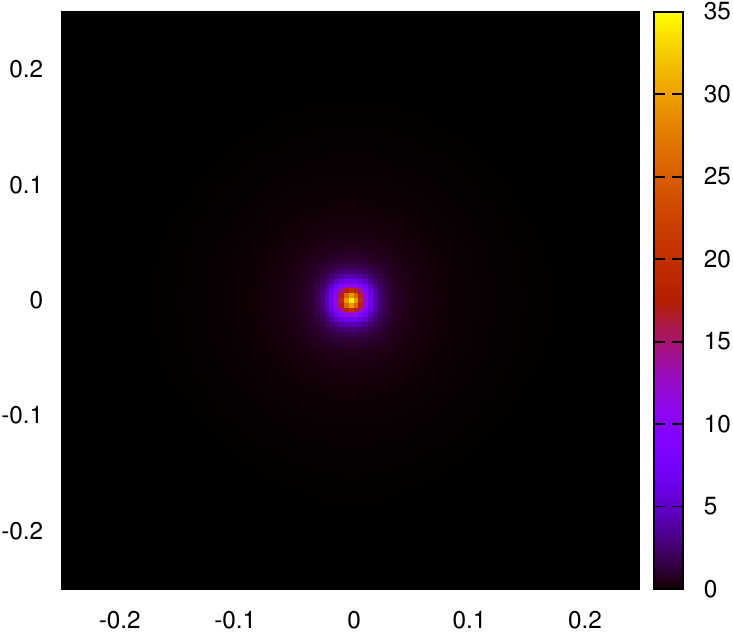} &
\includegraphics[width=4.5cm]{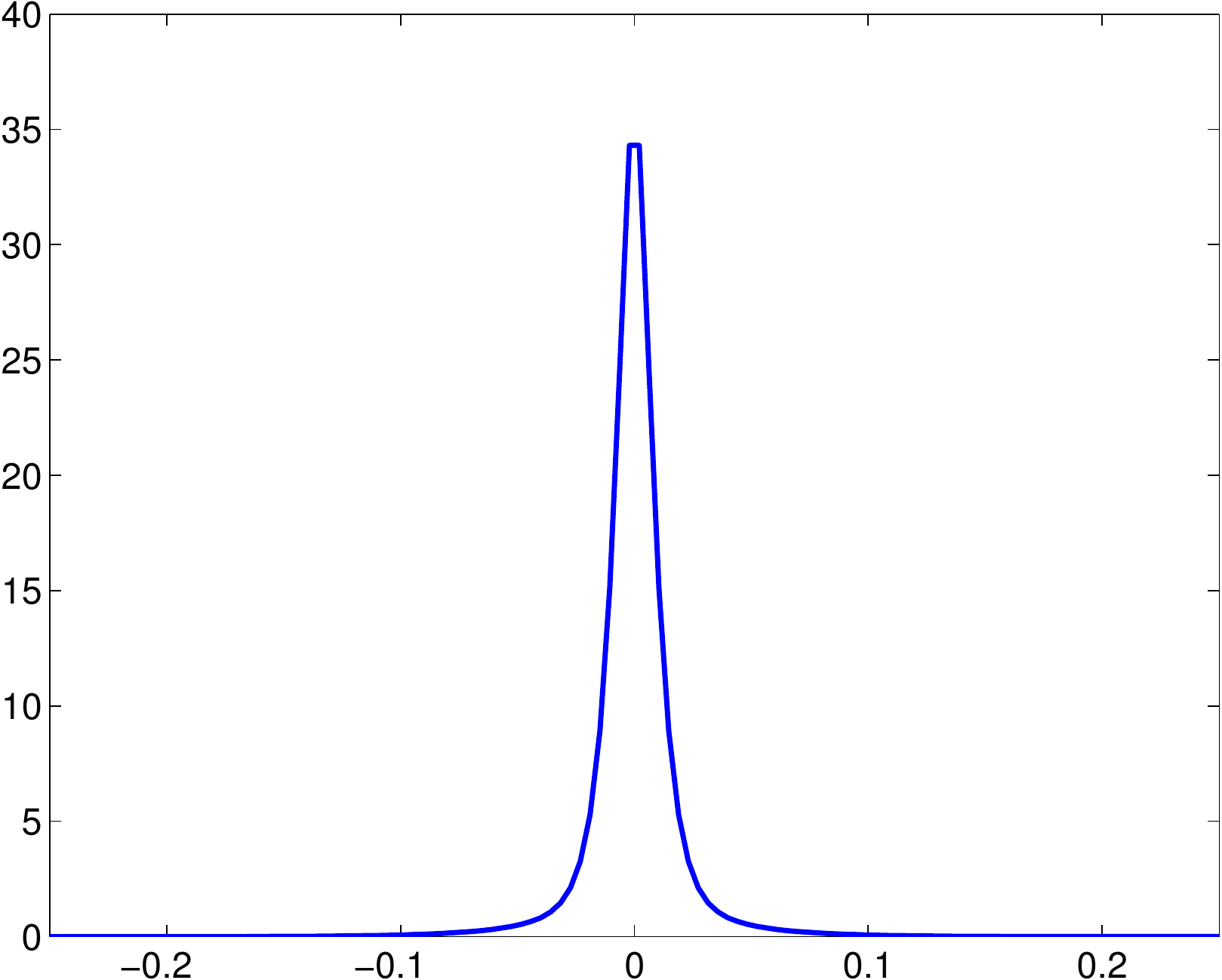} &
\includegraphics[width=4.5cm]{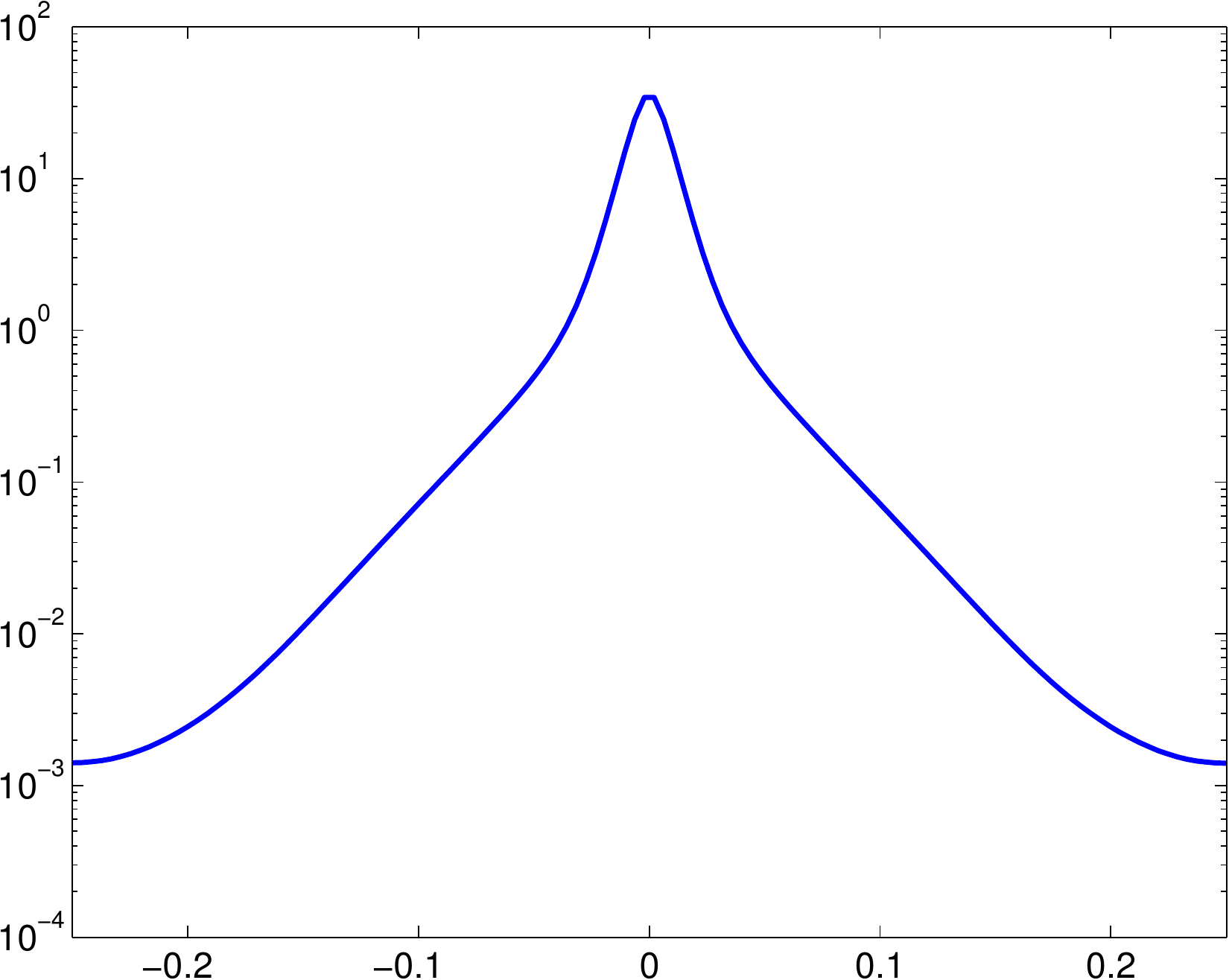} 
\\
$t=1\,\bar t$& $t=1\,\bar t$& $t=1\,\bar t$
\\
\includegraphics[width=4.2cm]{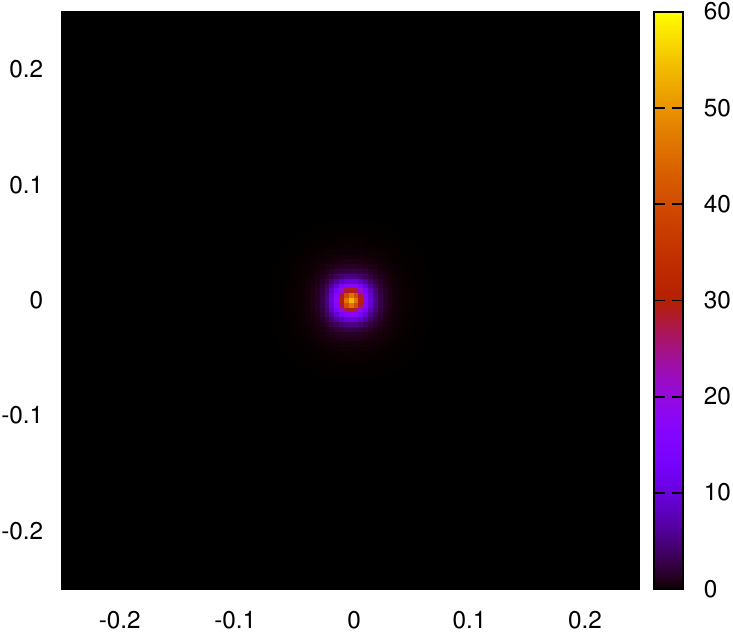} &
\includegraphics[width=4.5cm]{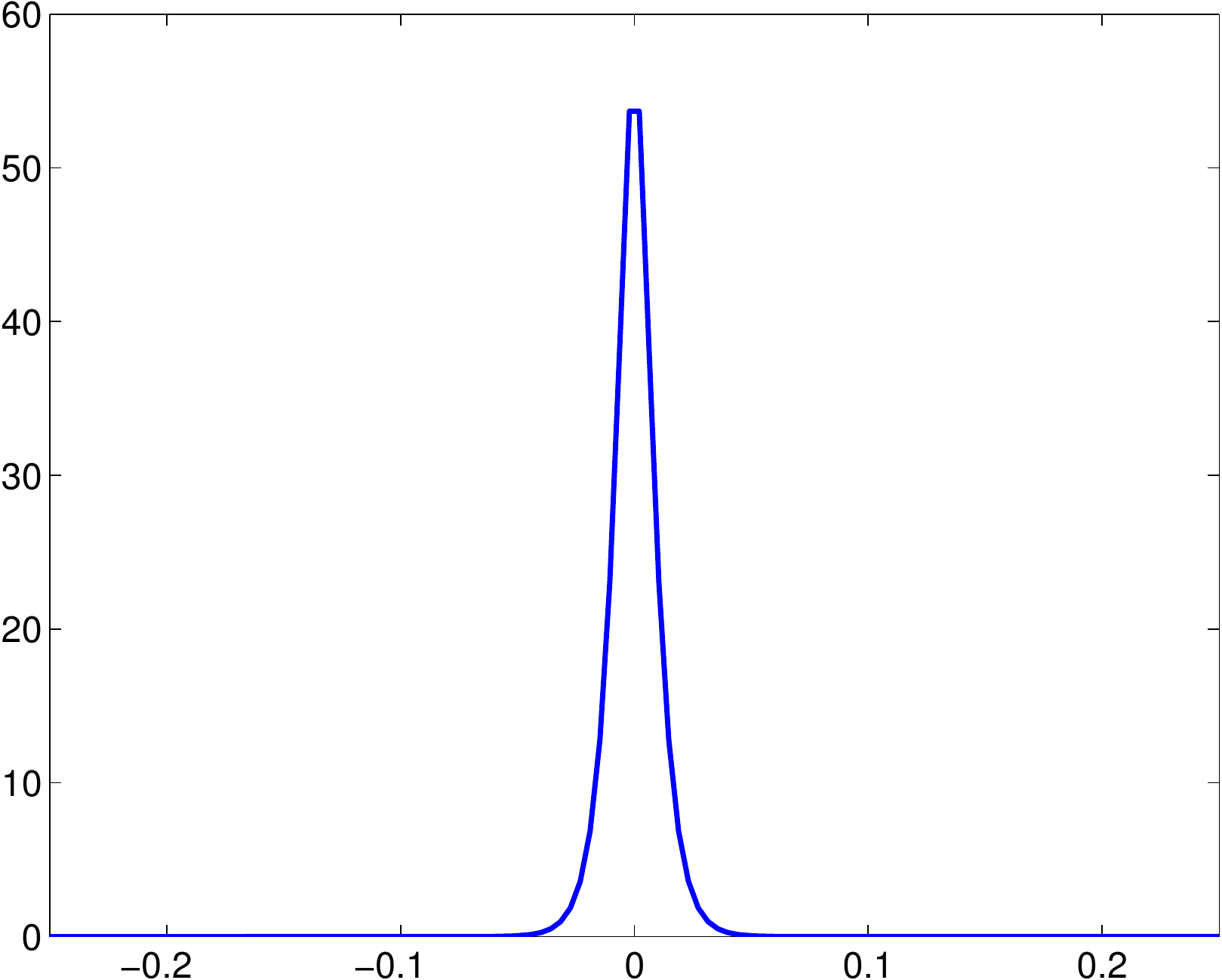} &
\includegraphics[width=4.5cm]{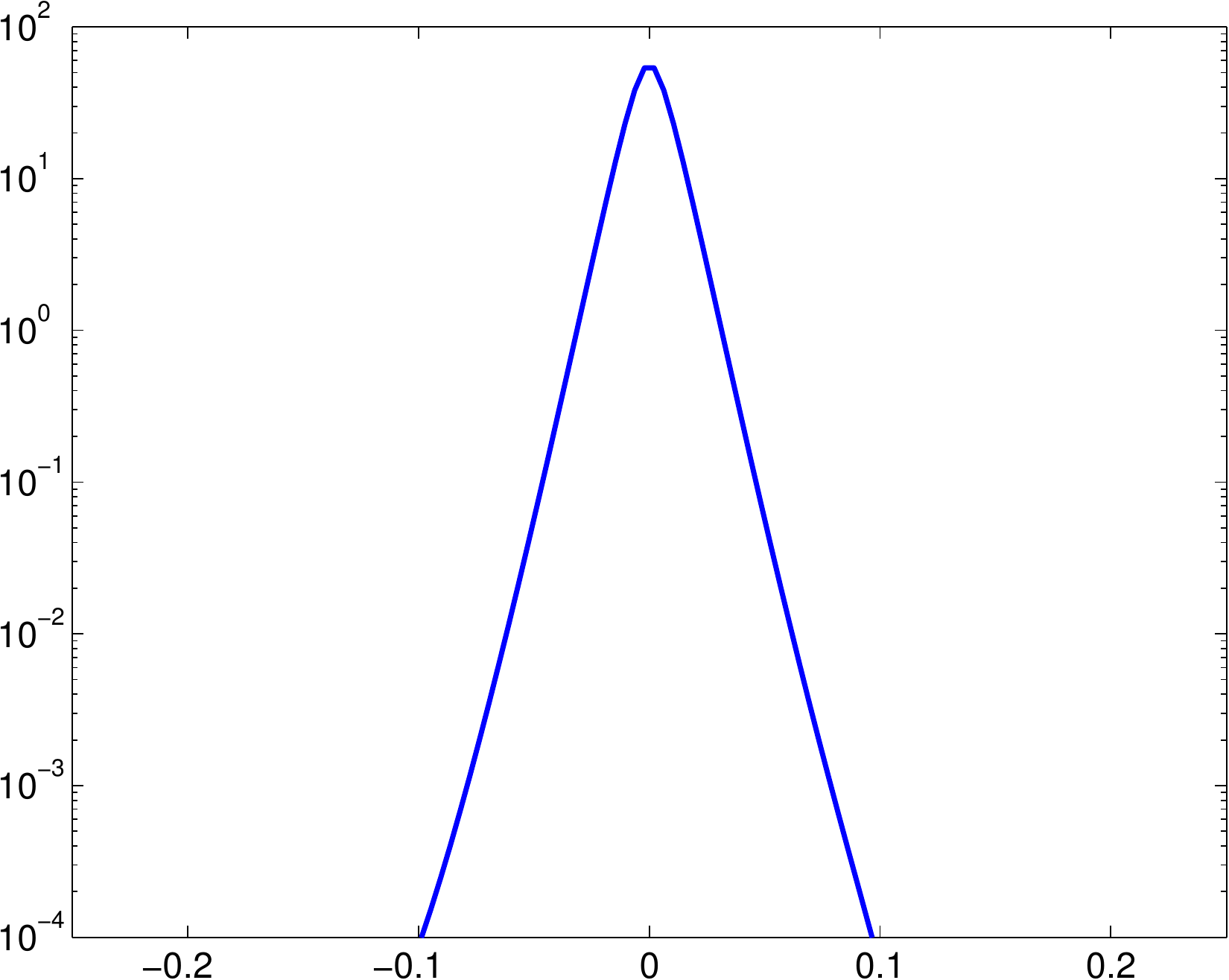} 
\\
$t=10\,\bar t$& $t=10\,\bar t$& $t=10\,\bar t$
\\
(a)&(b)&(c)
\end{tabular}
\caption{Test 1 : {\it Time evolution of the cell aggregation: (a) cell density in domain, (b) section plot of cell density, (c) section plot in log-scale of cell density. \label{Fig:test1-2}}}
\end{figure}

The evolution of the cell aggregation is presented in Figure~\ref{Fig:test1-2}. We observe that the initial density is a Gaussian distribution centered at $(0,0)$. At time $t=0.2\,\bar t$, the density forms a volcano profile, which disappears soon and is becoming into an exponential function at time  $t=1\,\bar t$. Finally at time $t=10\,\bar t$, a steady exponential function is formed. The last observation is similar with the sharp boundary of clusters in~\cite{bibMBBO}.

\begin{figure}[htbp]
\begin{tabular}{cc}
\includegraphics[width=7cm]{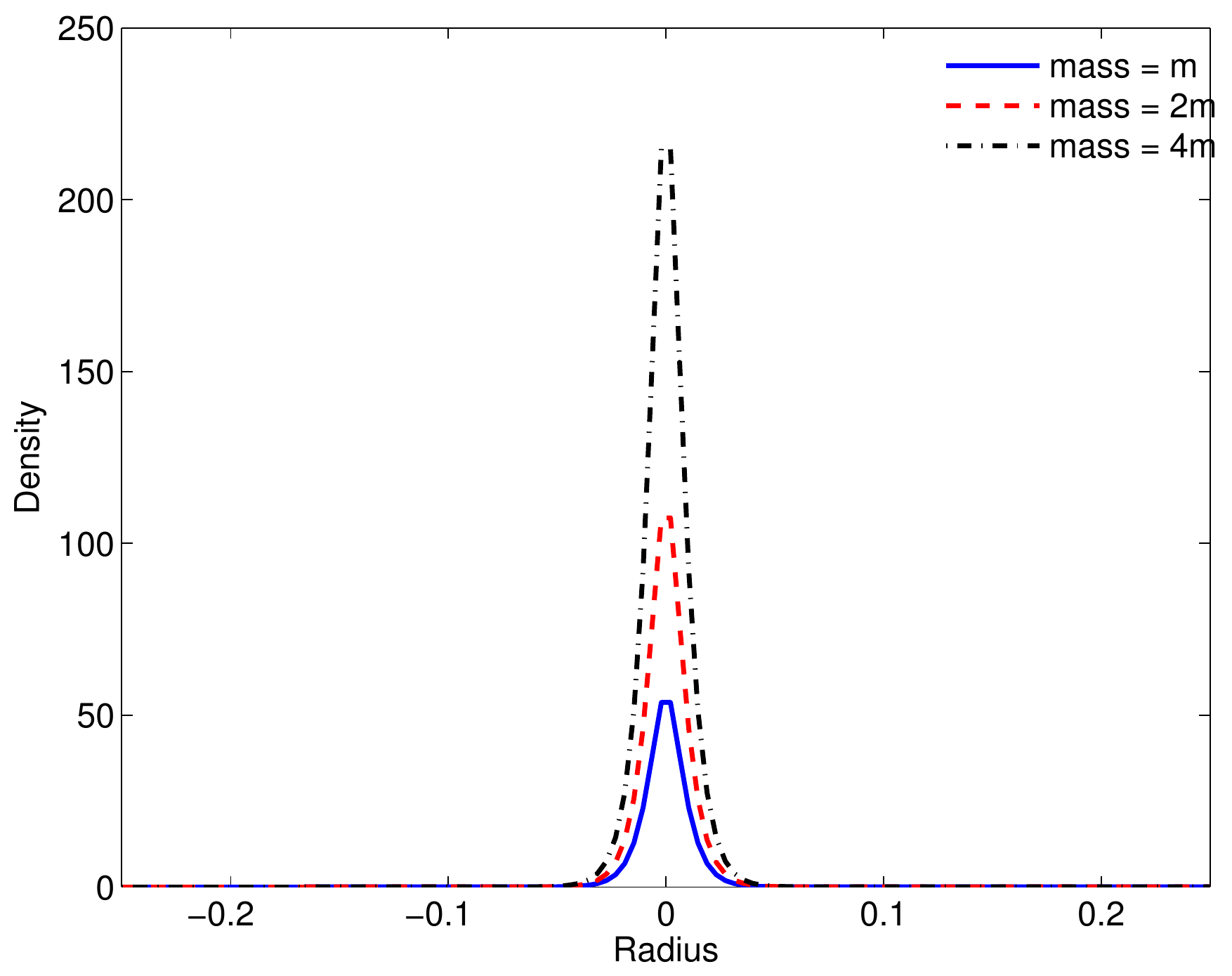} &
\includegraphics[width=7cm]{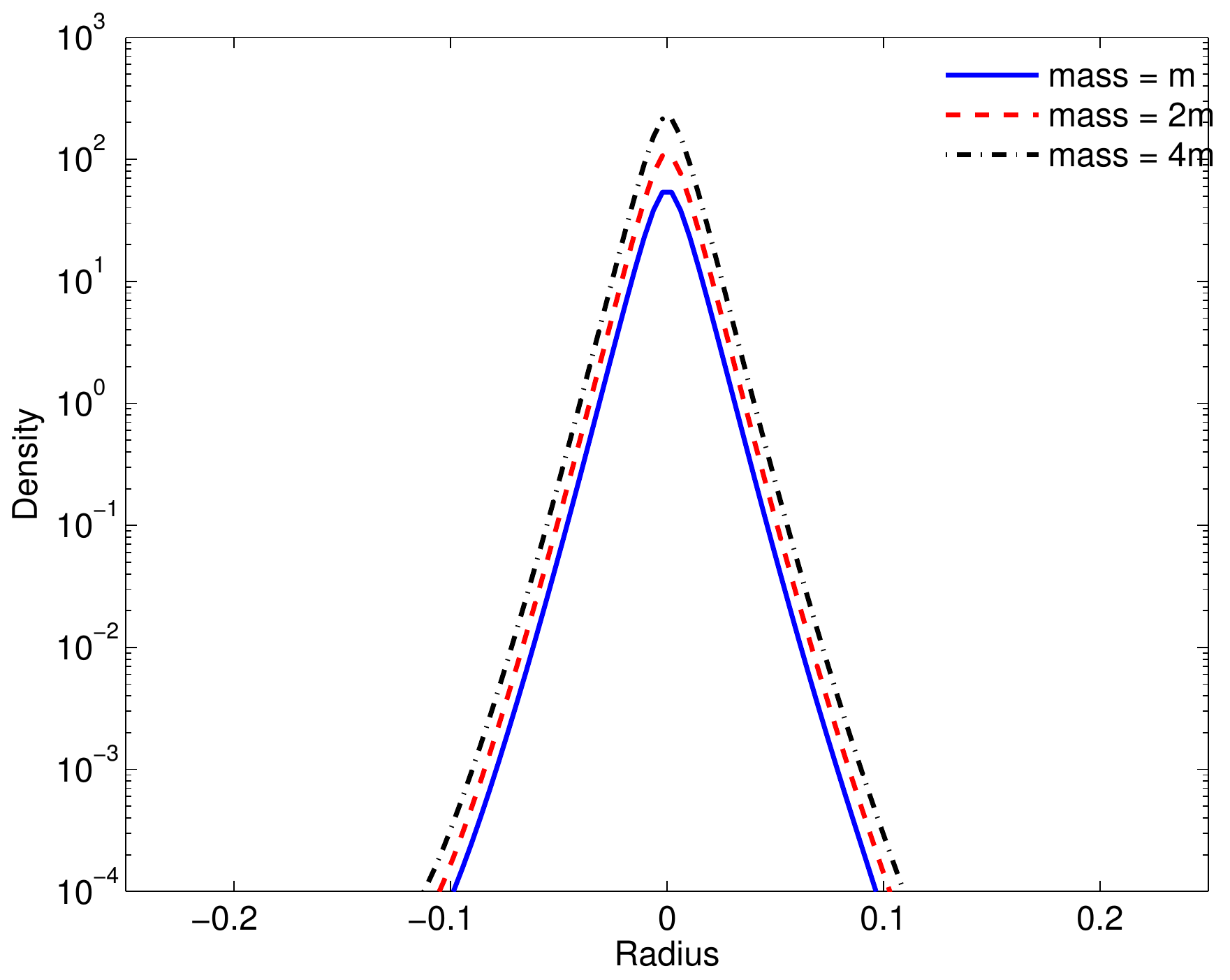}
\\
 (a)&(b)
\end{tabular}
\caption{Test 1 : {\it The steady-state size of clusters comparison for different total masses at time $t=10\,\bar t$:  (a) section plot of cell density, (b) section plot in log-scale of cell density. \label{Fig:test1-3}}}
\end{figure}

Furthermore, we observe in Figure~\ref{Fig:test1-3} that with the total mass equal to $m=\pi/100$, $\pi/50$ and $\pi/25$ respectively, the size of clusters is almost the same.
This phenomena coincides with the conclusion in~\cite{bibMBBO} that  the steady-state size of clusters is almost independent of the number of cells comprising them. In fact, following~\cite{bibSCBBSP} the steady density has a form
\begin{equation*}
\rho(\mathbf{x})\simeq \rho_0\exp(-\lambda |\mathbf{x}|), {\rm when }\, |\mathbf{x}|\rightarrow \infty,
\end{equation*}
where $\rho_0$ is the maximum density, $\lambda$ depends on $\chi_S$. Notice that the typical size of a cluster, defined as the mean radius given by 
\begin{eqnarray*}
  <|\mathbf{x}|>\,\,=\,\,\frac{\int_{\mathbf{x}}|\mathbf{x}|\rho(\mathbf{x})d\mathbf{x}}{\int_{\mathbf{x}}\rho(\mathbf{x})d\mathbf{x}} \,\,=\,\,\frac{1}{\lambda}
\end{eqnarray*}
is independent of the total mass. This computation shows that the size of a cluster only depends on the model parameters and is entirely independent of the initial mass condition.

\subsection{Test 2 :  wave propagation  in a disc}
Suspended Escherichia coli bacteria swim in convection-free geometries such as capillaries or micro-channels, collectively migrate towards nutrient-rich regions, in the form of propagation concentration waves~\cite{bibA}. In this test, we are interested in the model~\eqref{kinetic:eq}-\eqref{eq:BC:parabolique} to study concentration waves attracted by nutrient. Moreover, we consider a disc geometry  to verify the numerical discretization of boundary conditions in section~\ref{sec:discret_BC}.

The computational domain in space is  a square domain of size $[-3,3]^2$, which is uniformly divided by a mesh size of $n_x\times n_y=80\times80$. The disc is inside of the square with radius equal to $3$. It is clear that the boundary is not located on grids. Thus to achieve interior high order scheme, some artificial values on ghost points behind the boundary are needed. These values are given by using the numerical method presented in  section~\ref{sec:discret_BC}. Moreover the velocity space belongs to the unit circle  $S^1$, and is uniformly divided into $n_v=64$ parts. All the parameters are chosen as in~\cite{bibSCBPBS} and listed in Table~\ref{tab:parameters}.  The initial distribution function $f_0$ is given by a Gaussian function
\begin{equation*}
  f_0(\mathbf{x},\mathbf{v})=\rho_0\exp(|\mathbf{x}|^2),
\end{equation*}
where $\rho_0$ is constant. The initial chemoattractant $S_0$ is equal to 0 and the initial nutrient $N_0$ is homogeneous equal to 1.

\begin{figure}[htbp]
\begin{tabular}{ccc}
\includegraphics[width=5.cm]{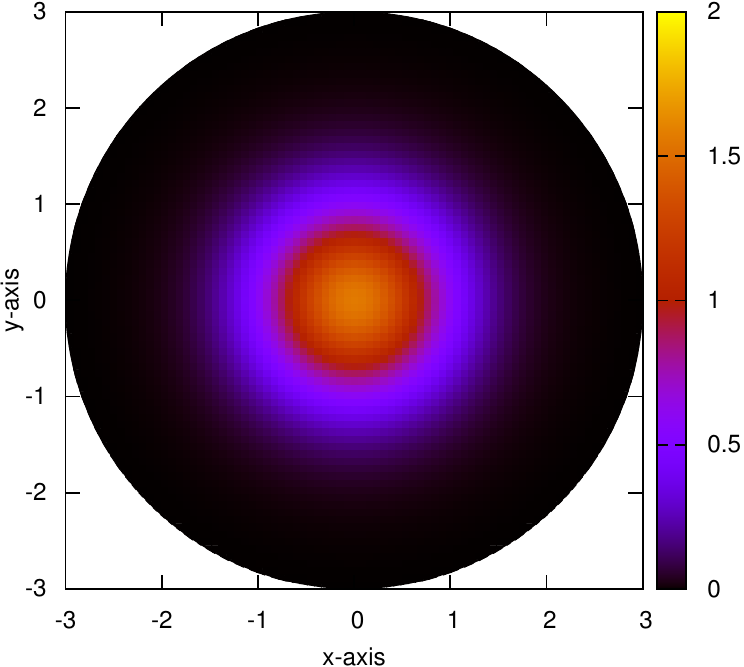} &
\includegraphics[width=5.cm]{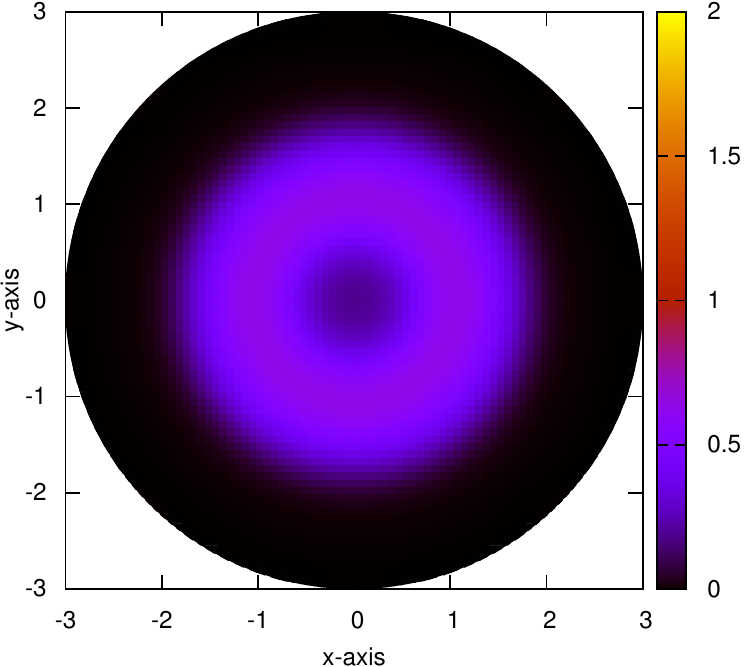} &
\includegraphics[width=5.cm]{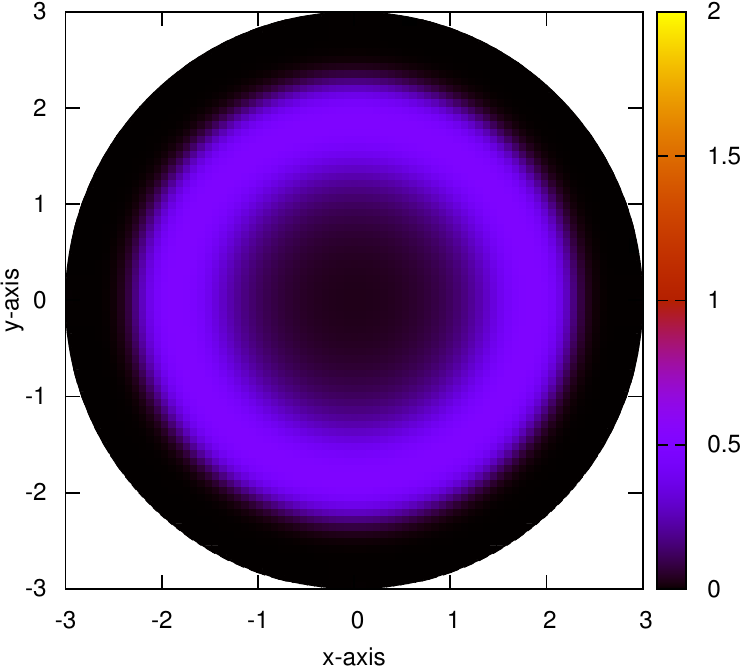} 
\\
$t=\,0$&$t=4\,\bar t$&$t=8\,\bar t$
\\
\includegraphics[width=5.cm]{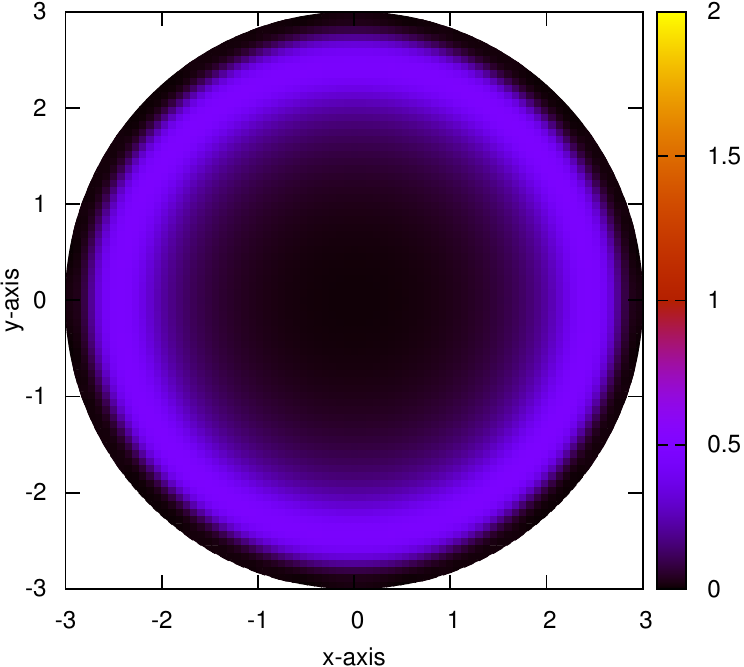} &
\includegraphics[width=5.cm]{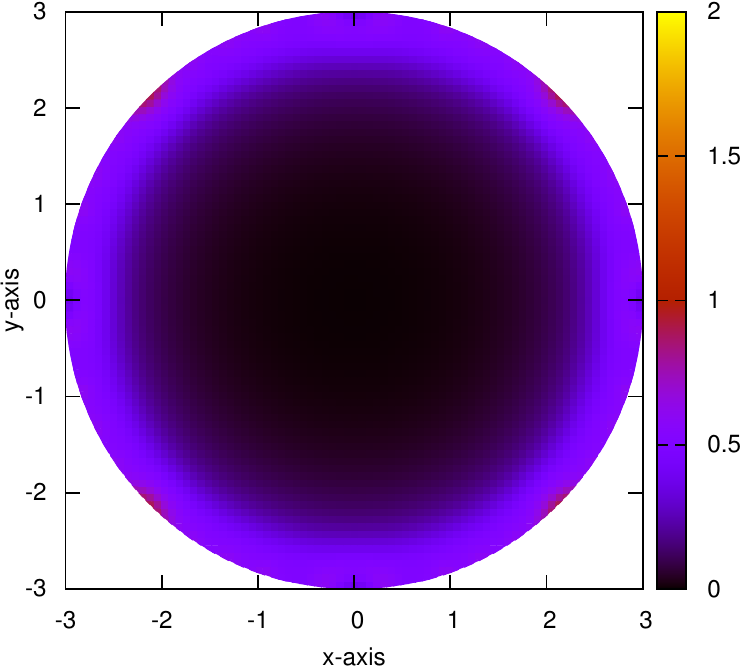} &
\includegraphics[width=5.cm]{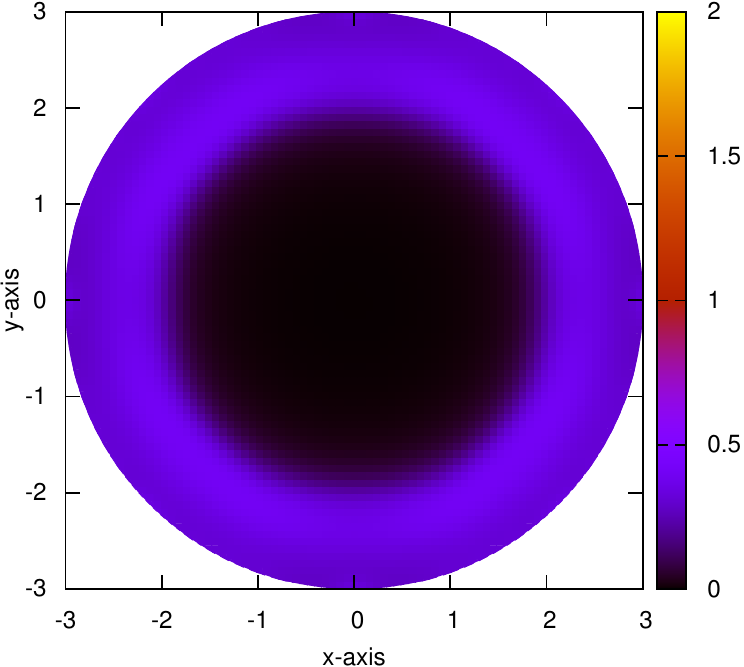} 
\\
$t=12\,\bar t$&$t=15\,\bar t$&$t=19\,\bar t$
\\
\includegraphics[width=5.cm]{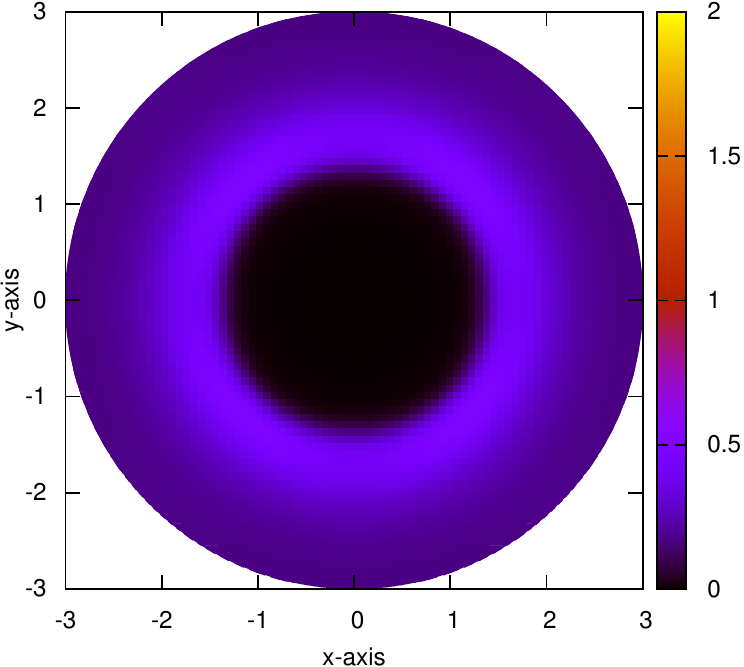} &
\includegraphics[width=5.cm]{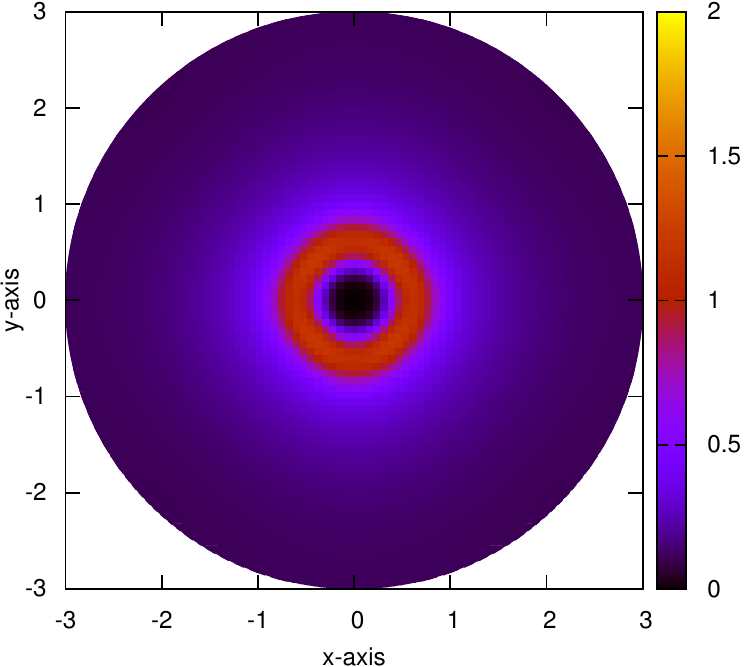} &
\includegraphics[width=5.cm]{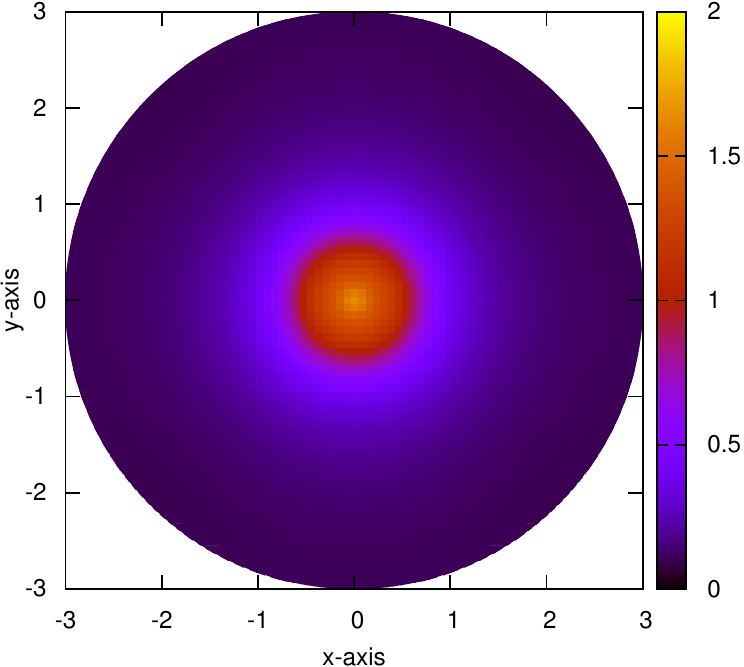} 
\\
$t=23\,\bar t$&$t=30\,\bar t$&$t=33\,\bar t$
\end{tabular}
\caption{\label{Fig2-1}Test 2 : {\it Time evolution of the cell density at different time.}}

\end{figure}

The time evolution of concentration wave in a disc is illustrated in Figure~\ref{Fig2-1}. In the first row of Figure~\ref{Fig2-1}, we observe that the initial Gaussian density is extending and  forming  a propagating wave. When the circle wave arrives at the boundary, all cell are reflected by the boundary and attracted by nutrient remained in the disc center (see the second row of  Figure~\ref{Fig2-1}). In the third row of Figure~\ref{Fig2-1}, the circle wave contracts back to disc center, and finally the cells concentrates at the disc center. We notice that cell diffusion appears when circle wave goes back to disc center. In fact, this diffusion is due to the stiffness of the response functions in the tumbling kernel~\cite{bibSCBPBS}.

\subsection{Test 3 : interaction of traveling waves in a $U$ shape.}
In this test, we focus on the influence of the reorientation on the shape. The simulations are compared to a particular experiment by courtesy of  Axel Buguin, Institut Curie (see the second row of Figure~\ref{Fig4-1}). We consider a channel of  $U$ shape, with initially homogeneous nutrient injected in the channel. Two clusters of bacteria are then imposed at two extremities of the channel. These two clusters move along the channel and finally meet at the channel center (the top of $U$ shape). We note that these two clusters keep their bar shape till their meeting.

To perform the simulation, we consider a channel of $U$ shape  with channel width  equal to 1 included in a rectangle computational domain $[0,8]\times[0,6]$, which is covered by an uniform mesh of size $n_x\times n_y=80\times60$. The velocity space belongs to  the unit circle $S^1$, and is uniformly divided into $n_v=64$ parts. The  numerical parameters used in the simulations are  listed in Table~\ref{tab:parameters}. The initial distribution function $f_0$ is given by a constant as follows
\begin{equation*}
  f_0(\mathbf{x},\mathbf{v})=\left\{
  \begin{array}{ll}
    0.25,&\text{if  } 0\leq y\leq1,\\[3mm]
    0,&\text{else}.
  \end{array}
    \right.
\end{equation*}
 The initial chemoattractant $S_0$ is equal to 0 and the initial nutrient $N_0$ is homogeneous equal to 1.

The numerical simulations are presented in the first row of Figure~\ref{Fig4-1}. In Figure~\ref{Fig4-1}(a), we observe that two bar shape clusters form and are moving to the channel center. In Figure~\ref{Fig4-1}(b), two clusters go forward along the half-ring channel and keep well their bar shape. Finally in Figure~\ref{Fig4-1}(c), two clusters meet at the channel center. We see that the numerical simulations have a good agreement with the experiment.

\begin{figure}[htbp]
\begin{tabular}{ccc}
\includegraphics[width=5.cm]{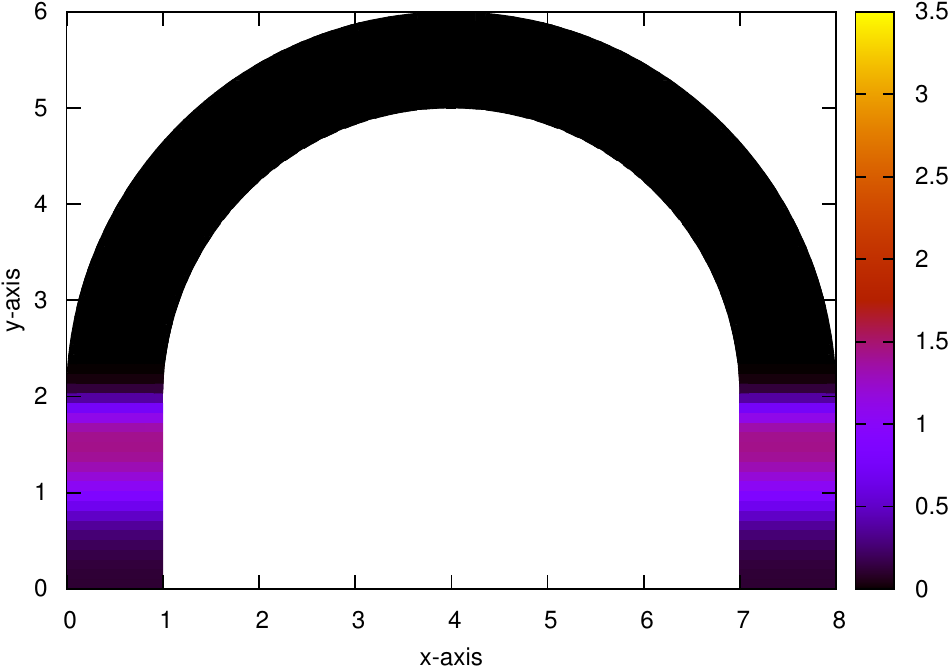} &
\includegraphics[width=5.cm]{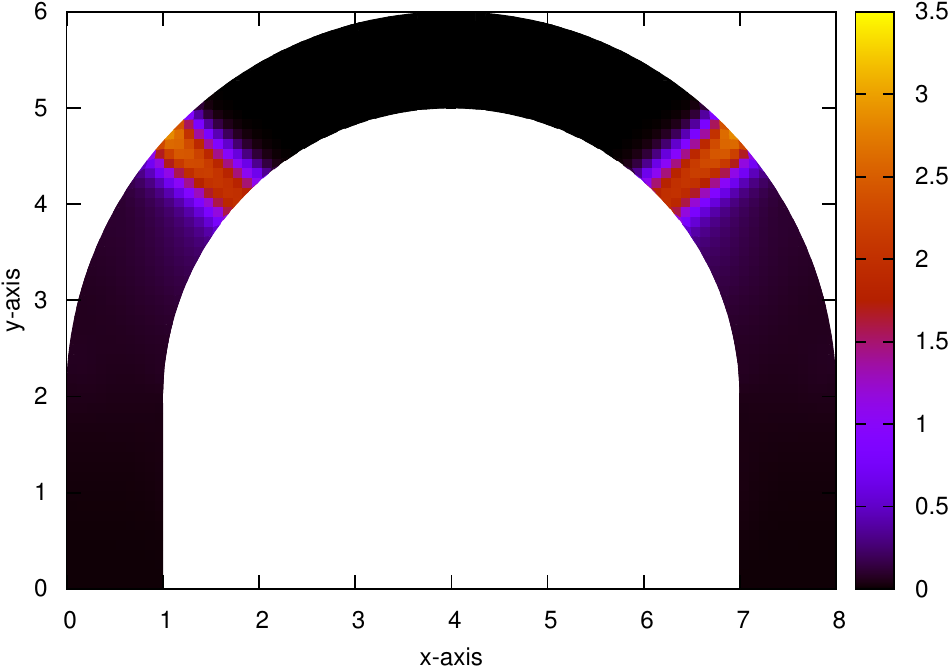} &
\includegraphics[width=5.cm]{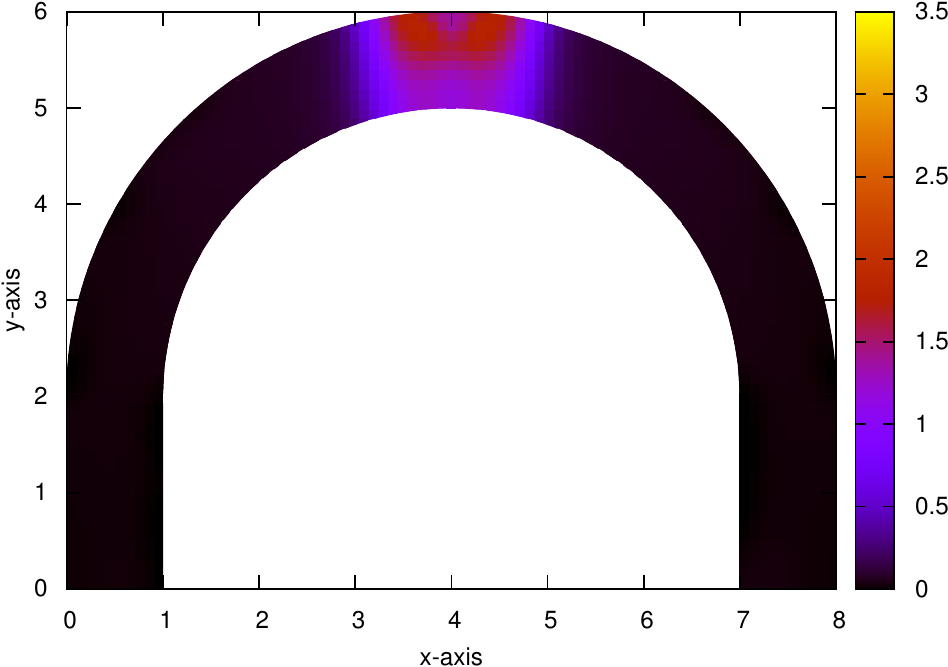} 
\\
\includegraphics[width=4.cm,height=4.cm]{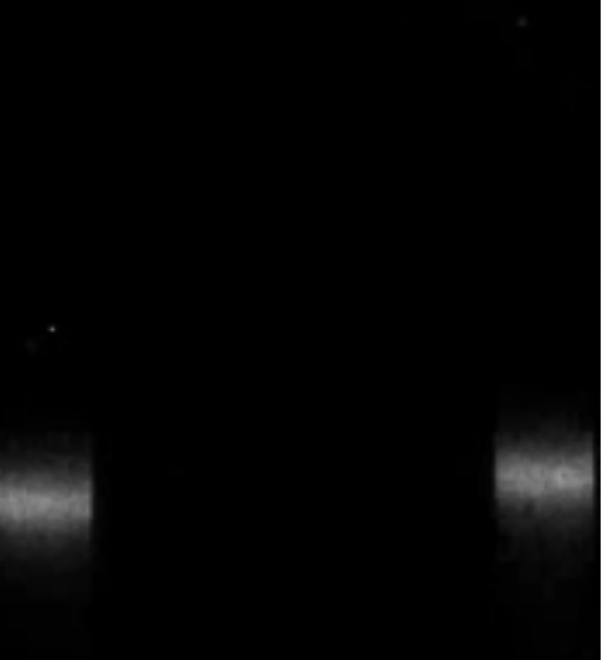} &
\includegraphics[width=4.cm,height=4.cm]{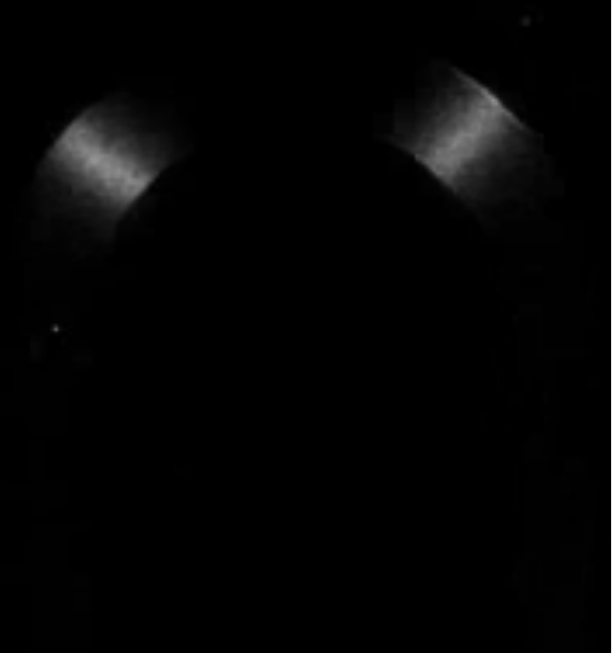} &
\includegraphics[width=4.cm,height=4.cm]{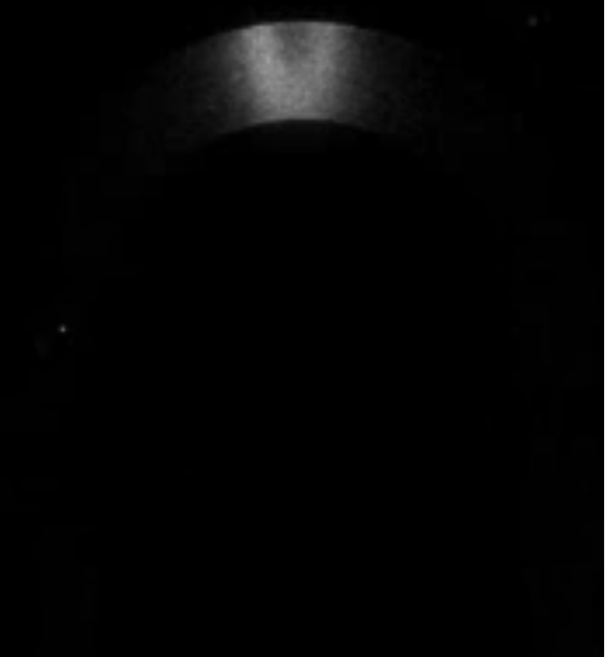} 
\\
(a)&(b)&(c)
\end{tabular}
\caption{Test 3 : {\it Time evolution of the cell density (a) t= 0.65 sec. (b) t = 2.65 sec. (c) t = 4.53 sec. At the top numerical results and at the bottom experiments on Escherichia coli (courtesy of Axel Buguin, Institut Curie).}}
\label{Fig4-1}
\end{figure}

\subsection{Test 4 : one traveling wave in a $U$ shape}
\begin{figure}[htbp]
\begin{tabular}{ccc}
\includegraphics[width=5.cm,height=6.cm]{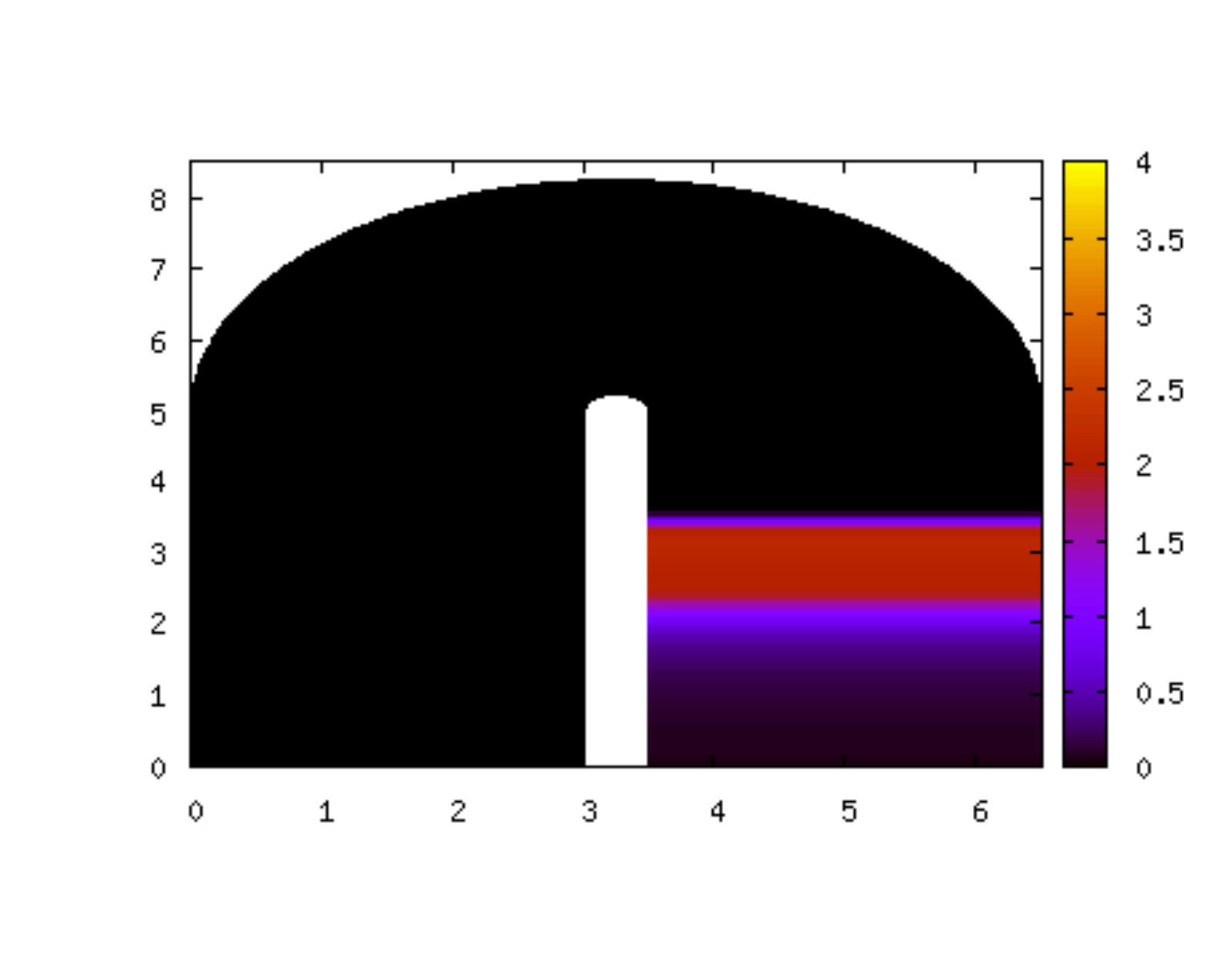} &
\includegraphics[width=5.cm,height=6.cm]{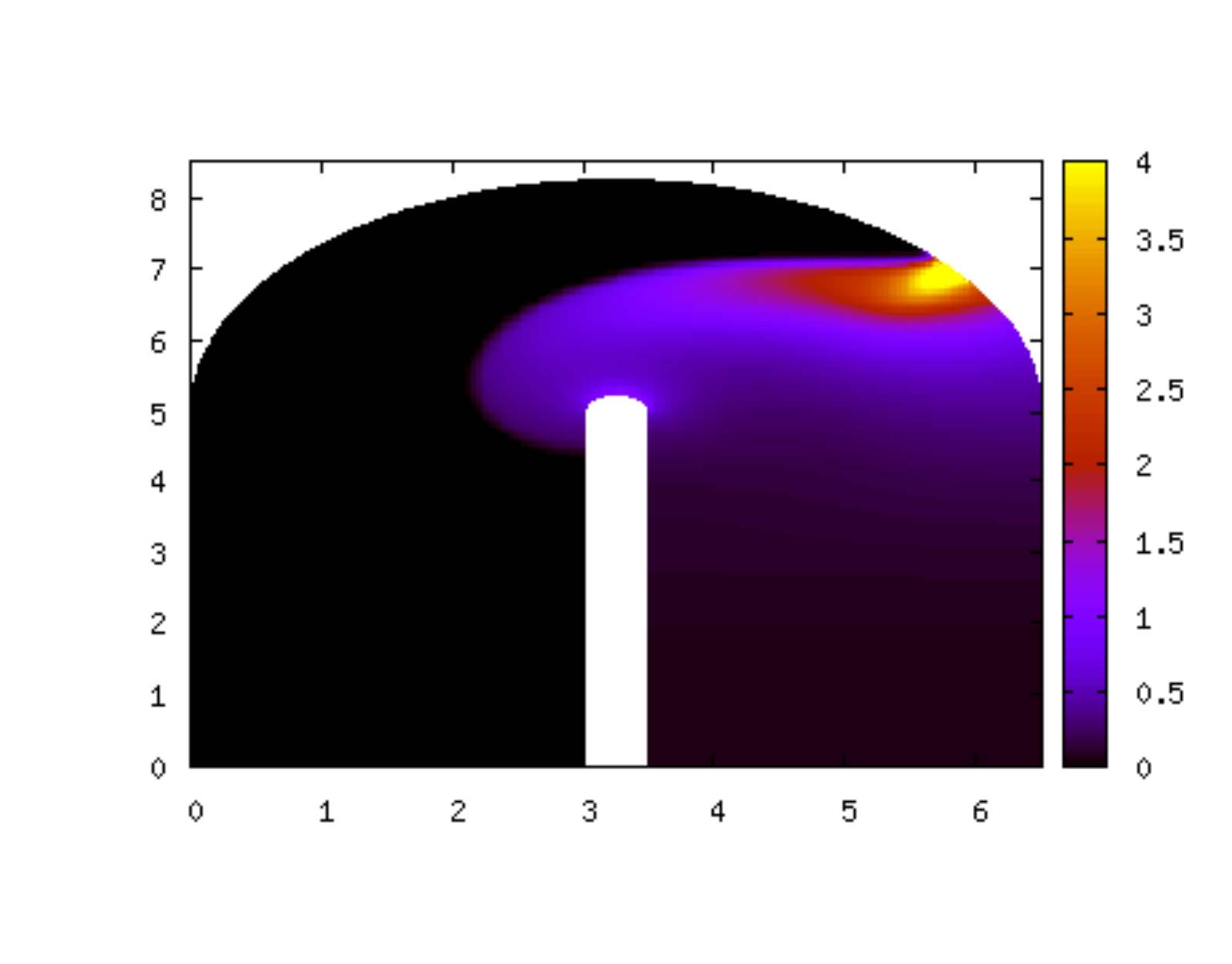} &
\includegraphics[width=5.cm,height=6.cm]{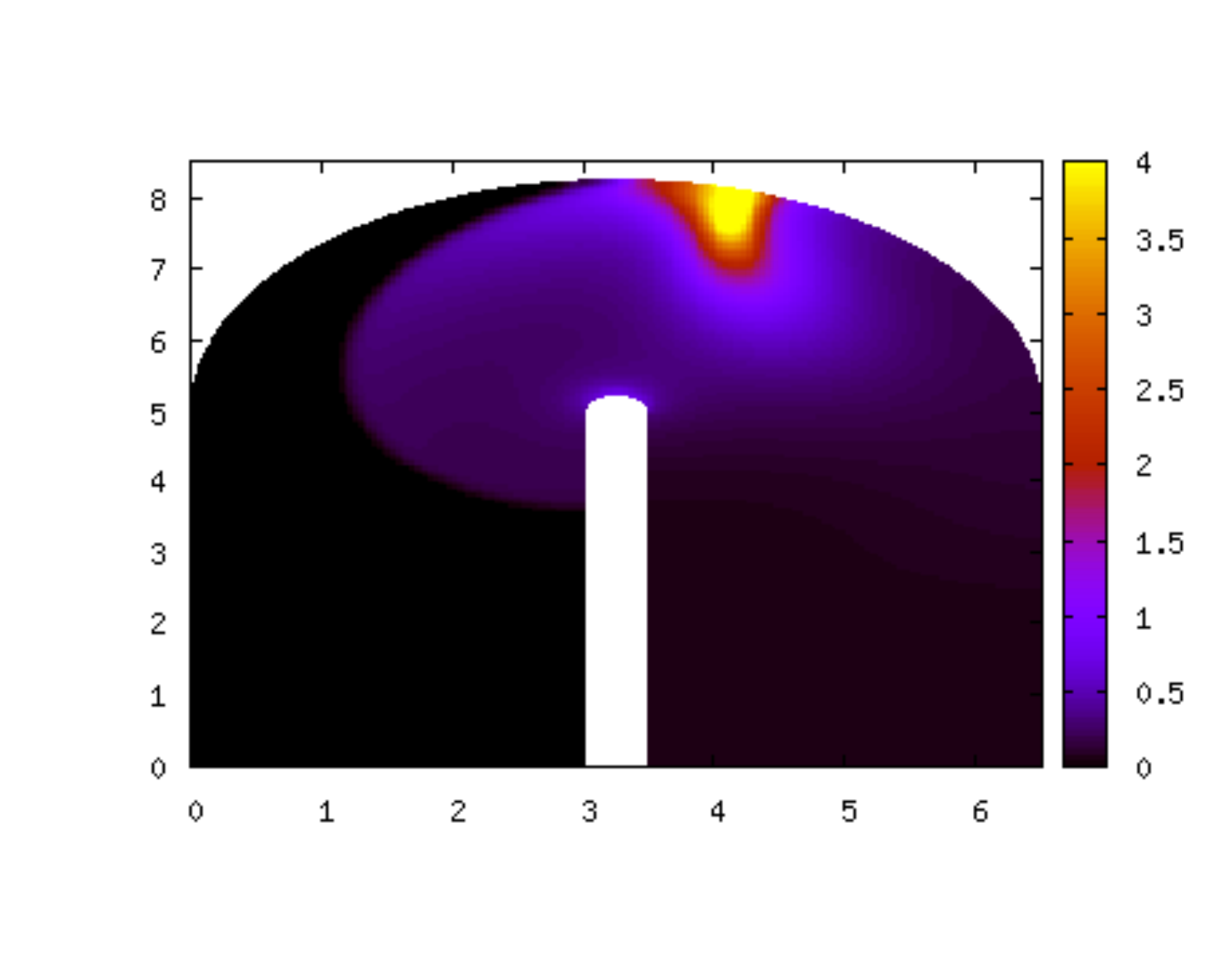} 
\\
\includegraphics[width=4.cm]{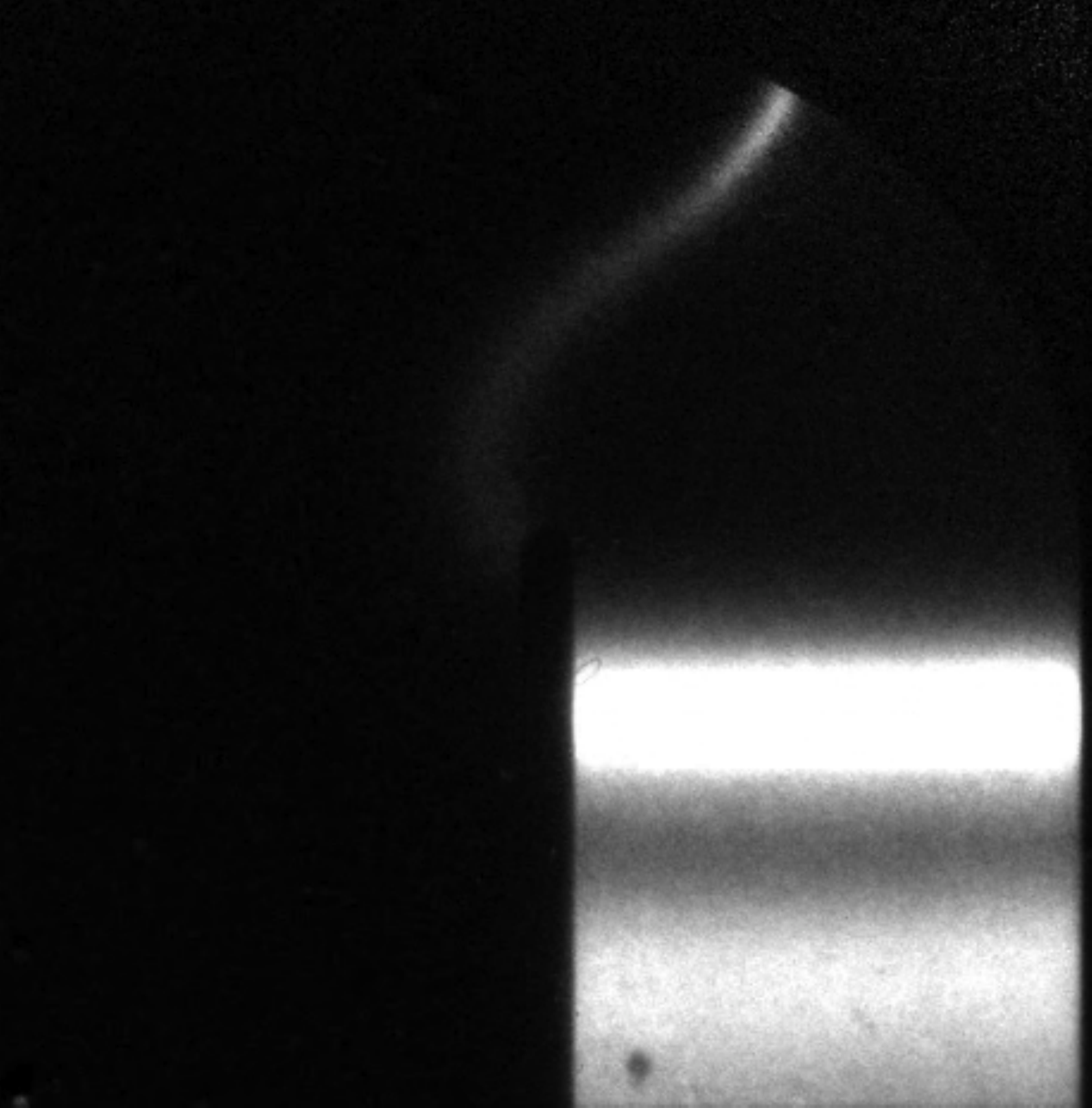} &
\includegraphics[width=4.cm]{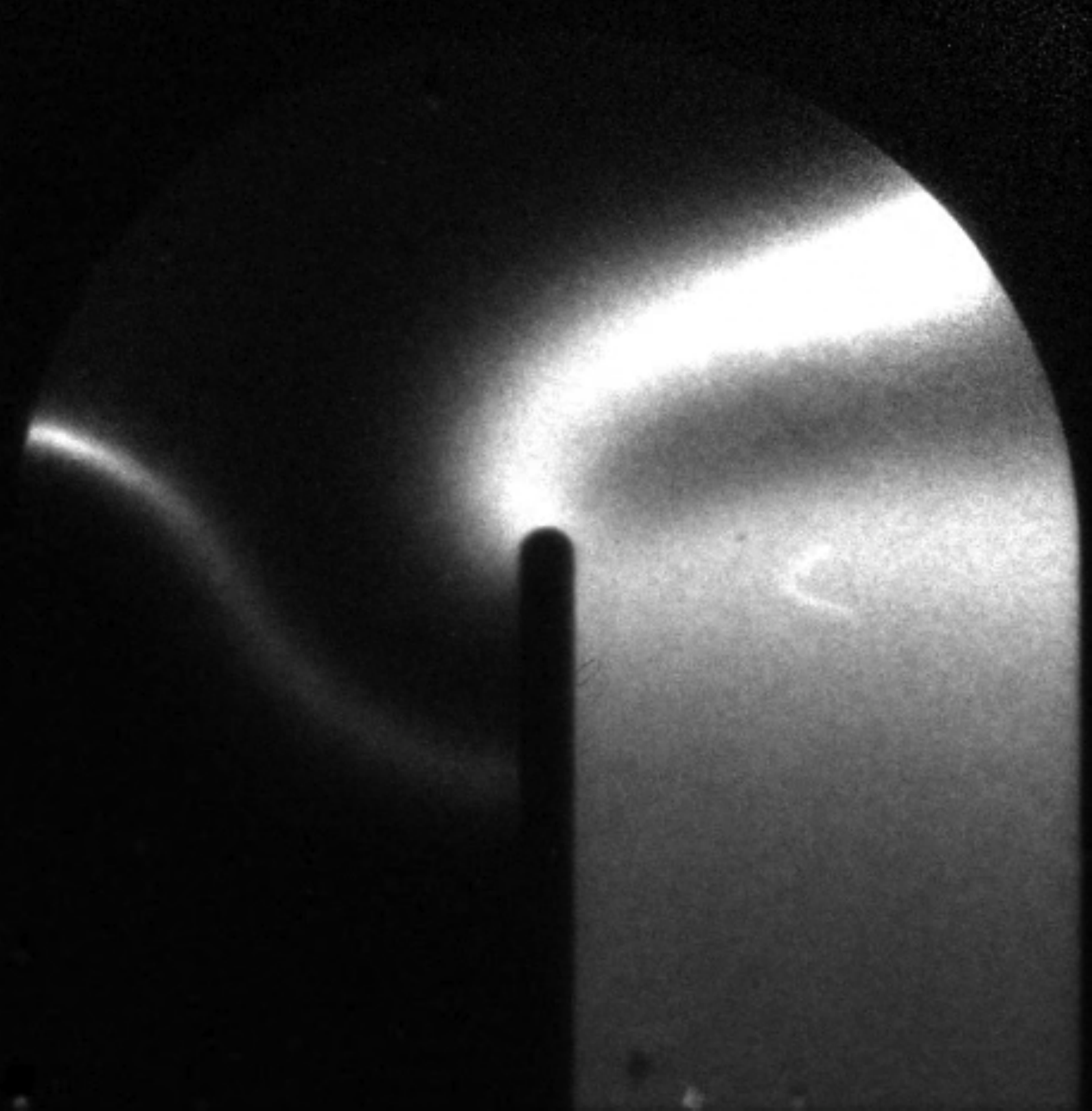} &
\includegraphics[width=4.cm]{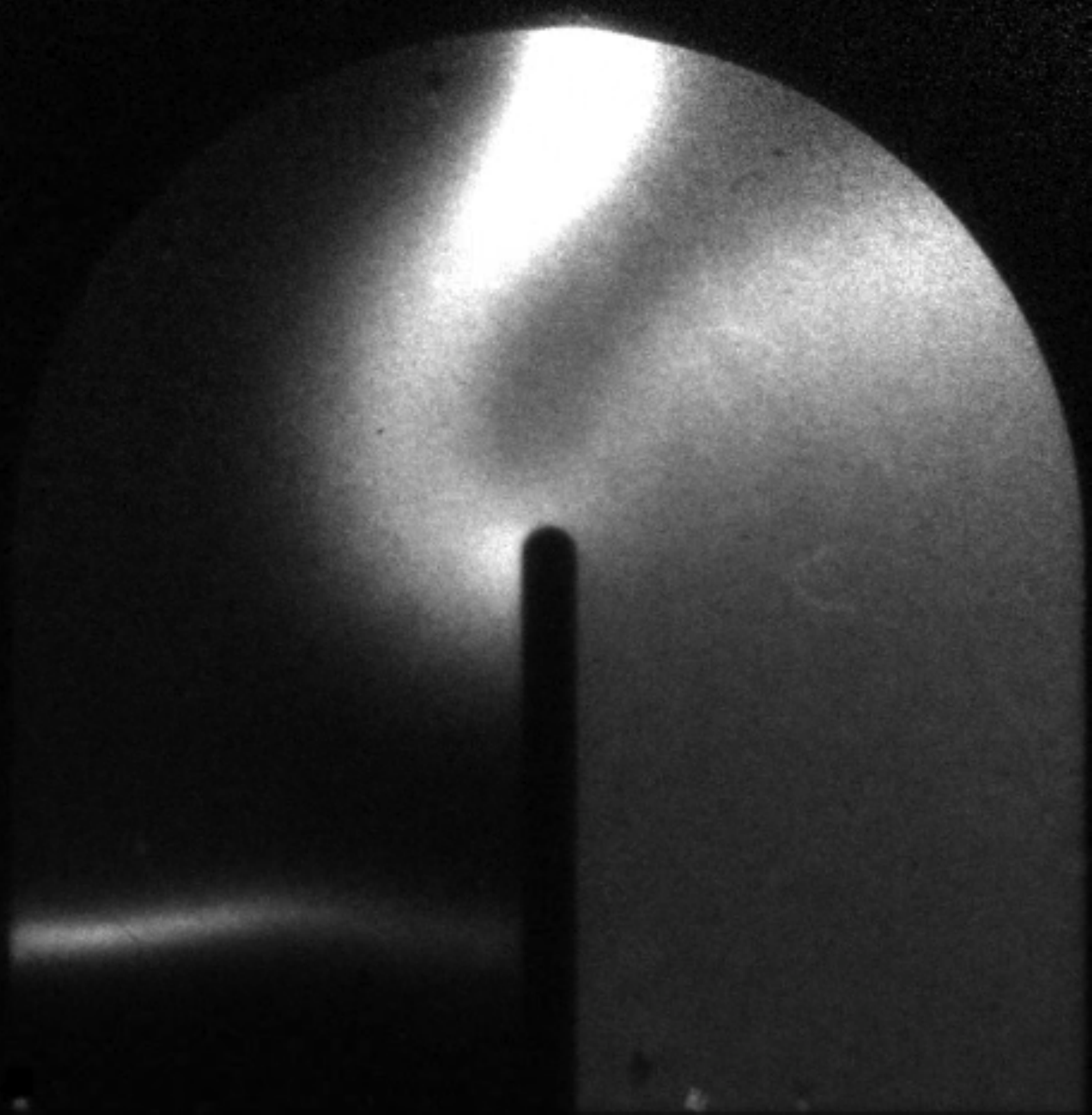} 
\\
(a)&(b)&(c)
\end{tabular}
\caption{Test 4 : {\it Numerical simulations (top) and experiments on Escherichia coli (bottom) : time evolution of the cell density (a) t = 7.5 $\overline{t}$ (b) t = 18.5 $\overline{t}$ (c) t = 23 $\overline{t}$ in sec. }}
\label{Fig5-1}
\end{figure}

\begin{figure}[htbp]
\begin{tabular}{ccc}
\includegraphics[width=5.cm,height=6.cm]{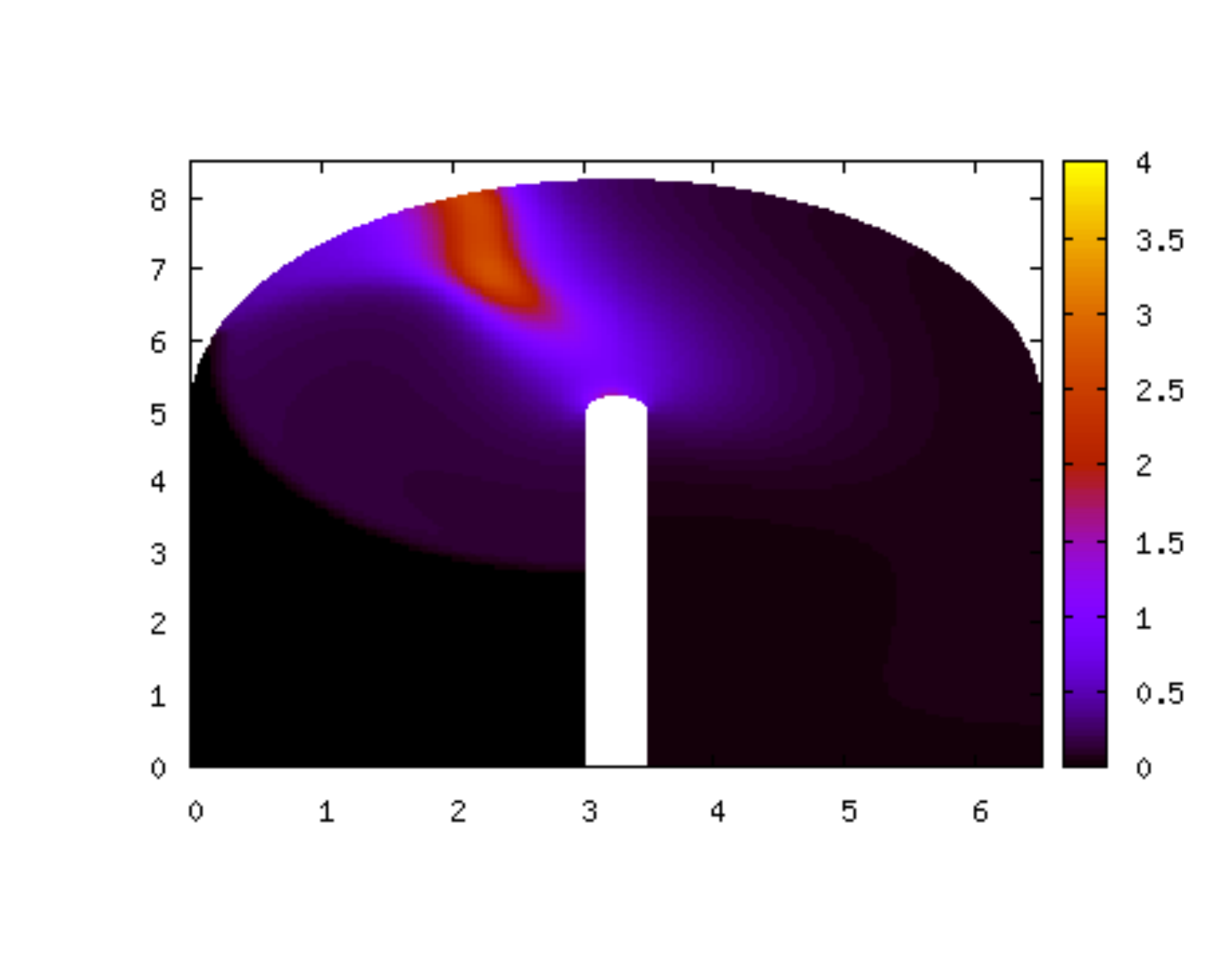} &
\includegraphics[width=5.cm,height=6.cm]{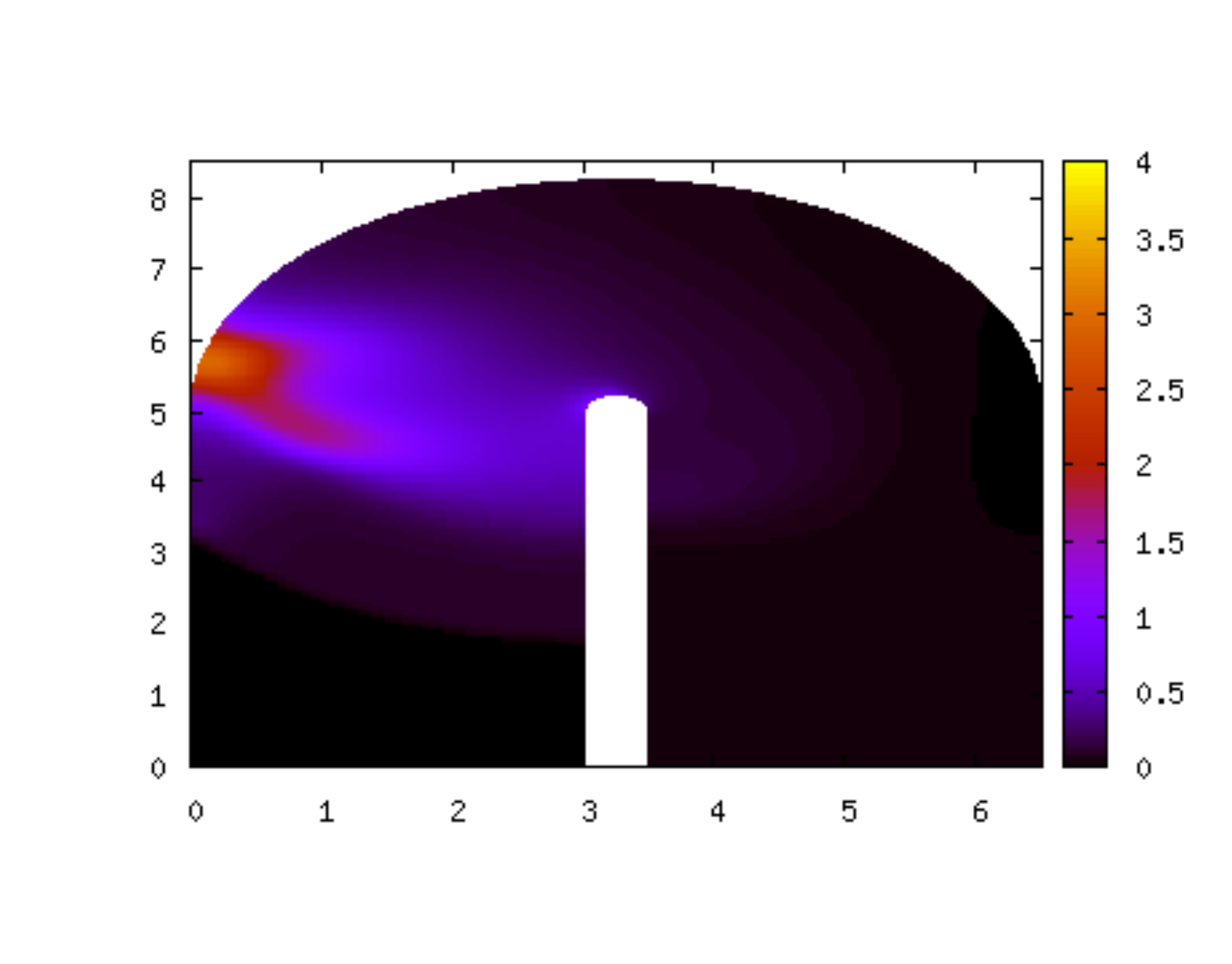} &
\includegraphics[width=5.cm,height=6.cm]{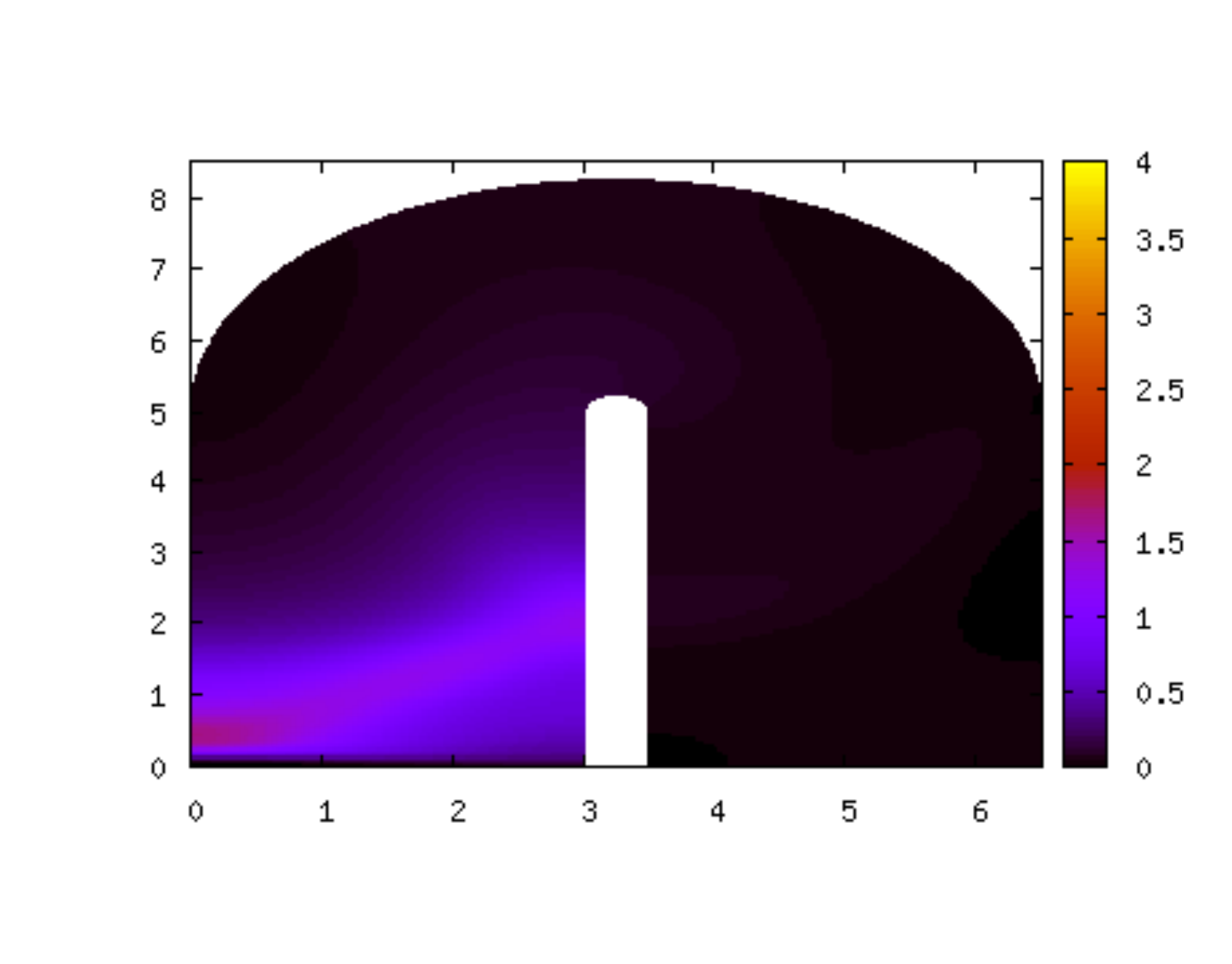} 
\\
\includegraphics[width=4.cm]{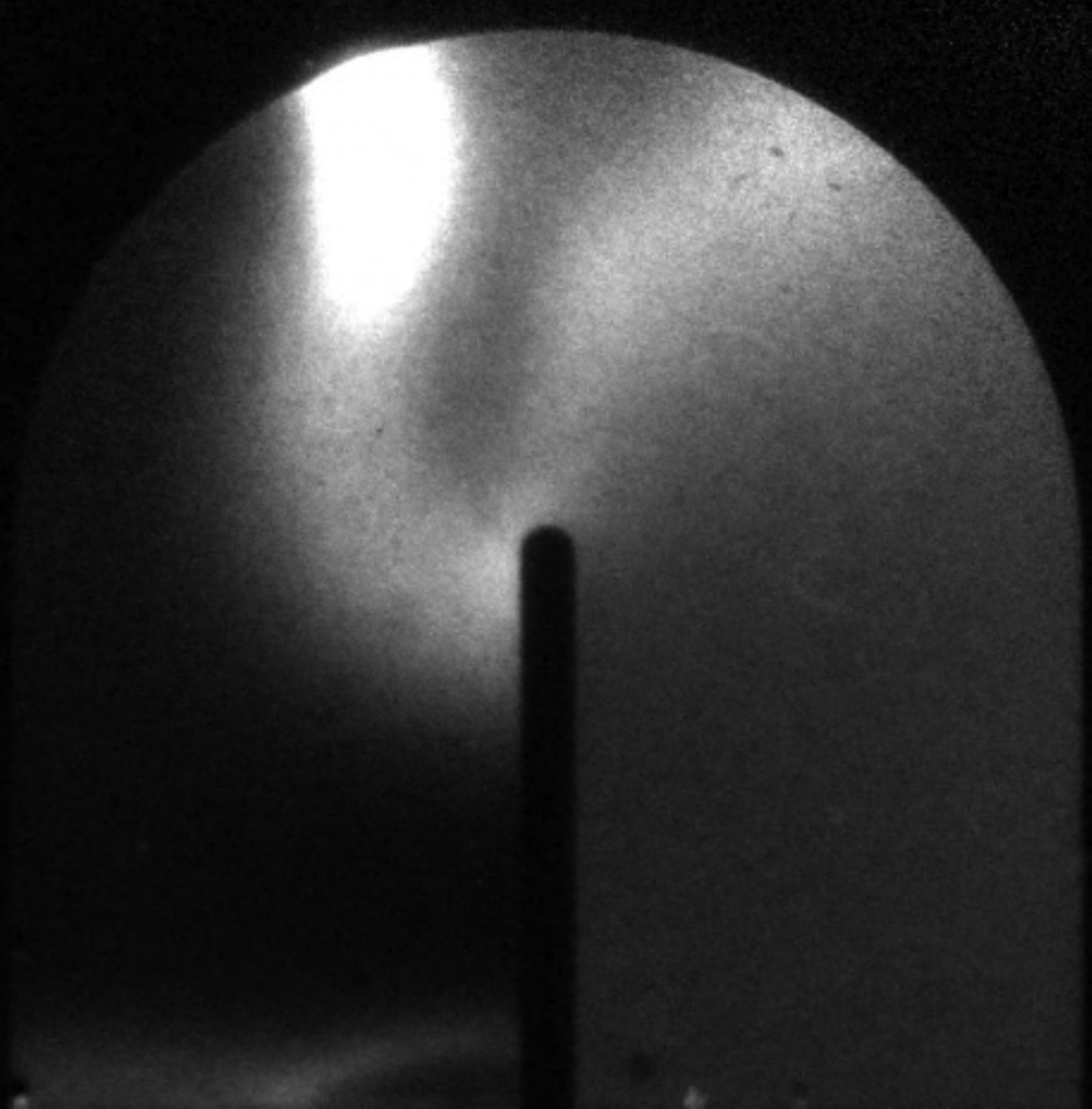} &
\includegraphics[width=4.cm]{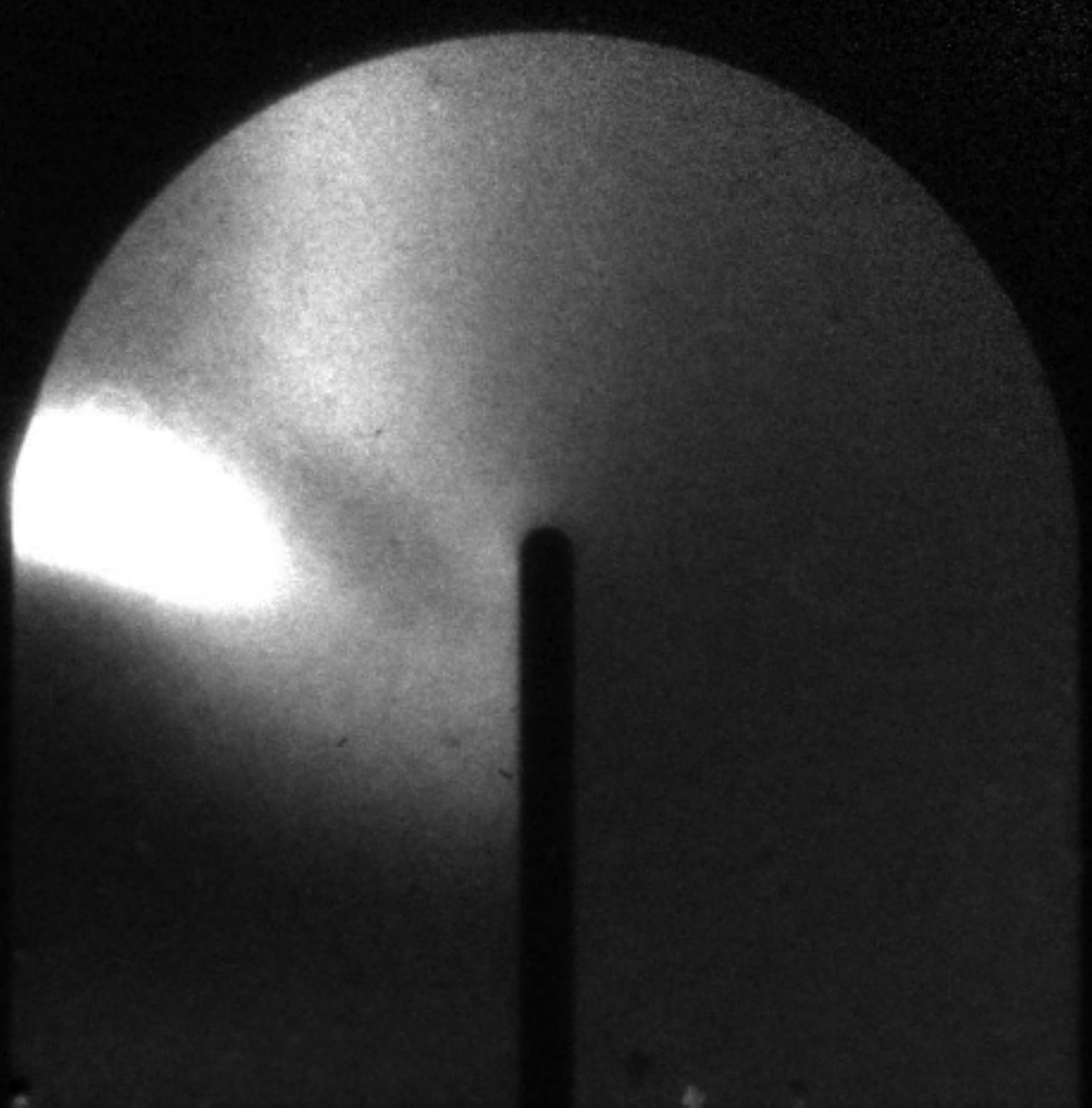} &
\includegraphics[width=4.cm]{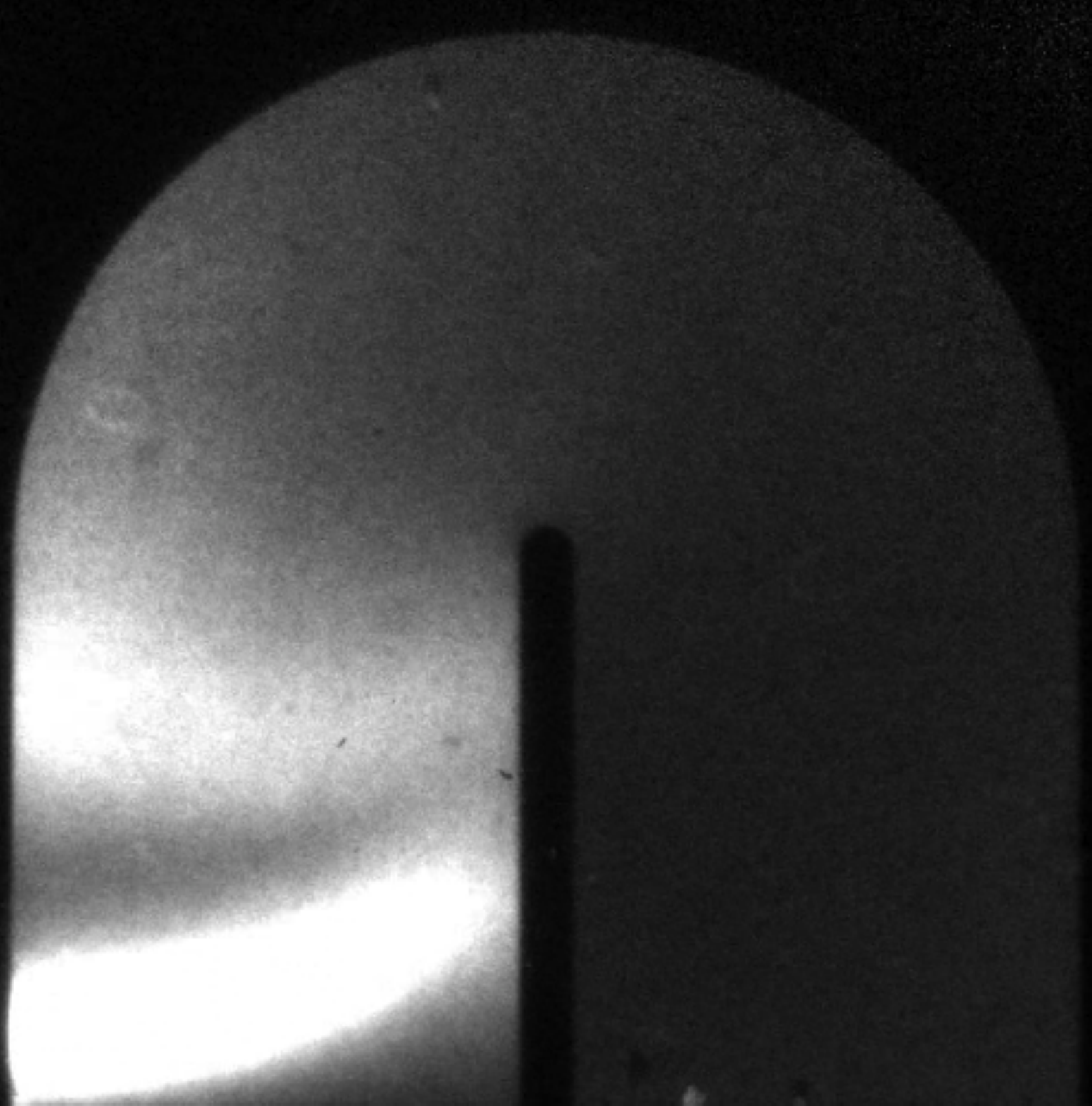} 
\\
(d)&(e)&(f)
\end{tabular}
\caption{Test 4 : {\it Numerical simulations (top) and experiments on Escherichia coli (bottom) :  time evolution of the cell density (d) t = 28 $\overline{t}$ (e) t = 34.5 $\overline{t}$ (f) t = 46.5 $\overline{t}$ in sec.}}
\label{Fig5-2}
\end{figure}

In the last test, we consider again traveling wave in a $U$ shape but with a  more wide channel than the previous one. The experience is shown in the second row of Figure~\ref{Fig5-1} and  Figure~\ref{Fig5-2}. Again we inject  homogeneous nutrient in the channel, but we consider only one cluster of bacteria at the right extremity of the channel. We observe that at straight part of channel the cluster goes ahead in a bar shape. Once the cluster enters into the half-ring part, the bacteria near interior circle goes faster than the one near the exterior circle. Moreover, bacteria contracts towards the exterior circle. Before the cluster enters into the left straight part of channel, bacteria concentrate almost near the exterior circle. When the cluster goes forward the other extremity, the cluster recovers its original bar shape.

To perform the simulation, we consider a channel of $U$ shape with channel width  equal to 3 included in a rectangle computational domain $[0,6.5]\times[0,8]$,  which is covered by an uniform mesh of size $n_x\times n_y=65\times80$. The velocity space belongs to  the unit circle $S^1$, and is uniformly divided into $n_v=64$ parts. The  numerical parameters used in the simulations are  listed in Table~\ref{tab:parameters}. The initial condition is the same as in the previous test.

The numerical simulations are presented in the first row of Figure~\ref{Fig5-1} and  Figure~\ref{Fig5-2}. We can see that the numerical simulations are very similar as the experience one. In fact, this phenomena is due to the directional persistence of chemotactic bacteria in a traveling concentration wave. When bacteria enter into a wide half-ring, they keep going straight and accumulate by the reflection of the exterior circle. It is very  similar like the effectiveness of centrifugal force, and can be observed significantly in a wide channel. This test shows that the model~\eqref{kinetic:eq}-\eqref{eq:BC:parabolique} represents well the chemotactic bacteria behavior and our numerical discretization based on Cartesian mesh approximates well the continuous model.

\subsection{Test 5 : long time behavior and pattern formations.}

In this last numerical test, we consider the full kinetic model with cell divisions and degradation
\begin{equation}
\label{kinetic:eq1}
  \frac{\partial f}{\partial t} \,+\, \mathbf{v}\cdot\nabla_\mathbf{x} f \,=\,Q(f)\,+\,r f - \gamma\Theta(\rho)\,\Theta(\rho-\rho_{\infty}) , \quad  \mathbf{x}\in \Omega,\,  \quad \mathbf{v}\in V,
\end{equation} 
where the division rate is given by \cite{PPLB}
$$
r = \frac{G_0\,N}{\sigma+N}, 
$$
it takes into account two experimental facts : the slowing down of the growth rate for low nutrients concentrations and its finite quantity  for high concentrations. The third term on the right hand side of (\ref{kinetic:eq1}) describes the vegetative cells into anabiotic form due to the increase of the local density. The transition starts when the total cell density reaches the value $\rho_\infty$ and $\Theta$ is the Heaviside function.

Actually such a source term has been proposed in \cite{PPLB}   in the framework of a macroscopic Patlak-Keller-Segel model where the unknown is the total density $\rho$. Numerical simulations of this macroscopic model have shown pattern formations as the ones observed in experiments \cite{PPLB,Berg}.  Here we consider the following initial density
$$
f_0(\mathbf{x},\mathbf{v}) \,:=\, 
\left\{
\begin{array}{ll}
1 & {\rm if}\, |\mathbf{x}|\leq  1,
\\
0 & {\rm else}\,
\end{array}\right.  
$$
and the initial nutrient concentration is uniform $N_0=0.5$. In the nutrient and chemoattractant system of equations (\ref{NS:eq}) we take the following parameters $D_S=D_N=1$, the production rate of chemoattractant $b=20$, the degradation rate of chemoattractant $a=8$ and the consumption rate of nutrient $c=0.8$, whereas in (\ref{kinetic:eq1}), we choose $\sigma=0.1$, $\rho_\infty=15$ and for the turning operator (\ref{def:T})-(\ref{def:psi}),we have  $\delta=20$, $\psi_0=1$, $\chi_N=1/2$ and $\chi_S=1/10$. The numerical simulations are presented in Fig.~\ref{Fig4-1}. We observe that for high initial nutrient concentration, cell density in the expanding ring becomes sufficient both for the break of stability of the uniform cell distribution and for their aggregation. In particular, if after the formation of the successive set of aggregates the expanding ring had time to grasp a certain part of the cells, participating in aggregation, then a radial pattern is formed.
\begin{figure}[htbp]
\begin{tabular}{ccc}
\includegraphics[width=5.cm,height=5.cm]{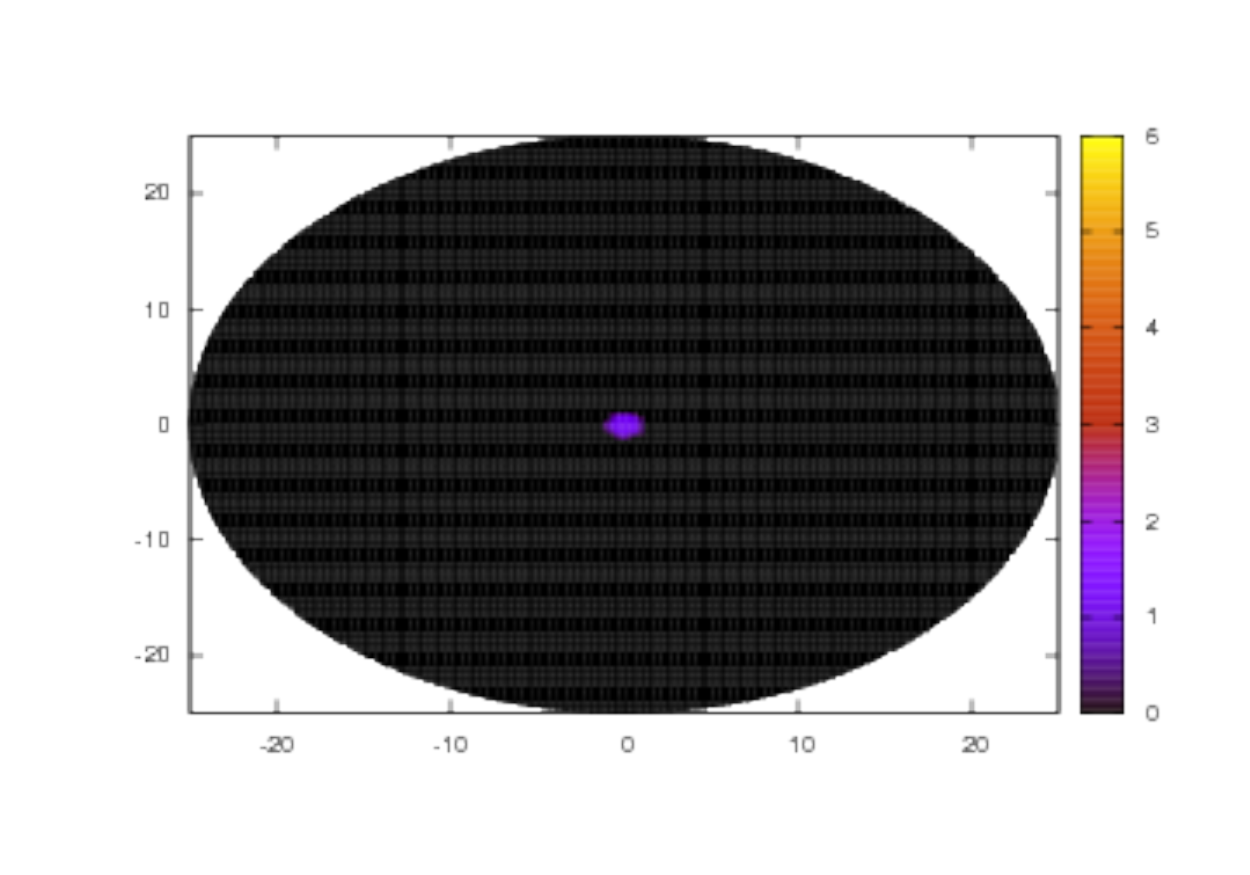} &
\includegraphics[width=5.cm,height=5.cm]{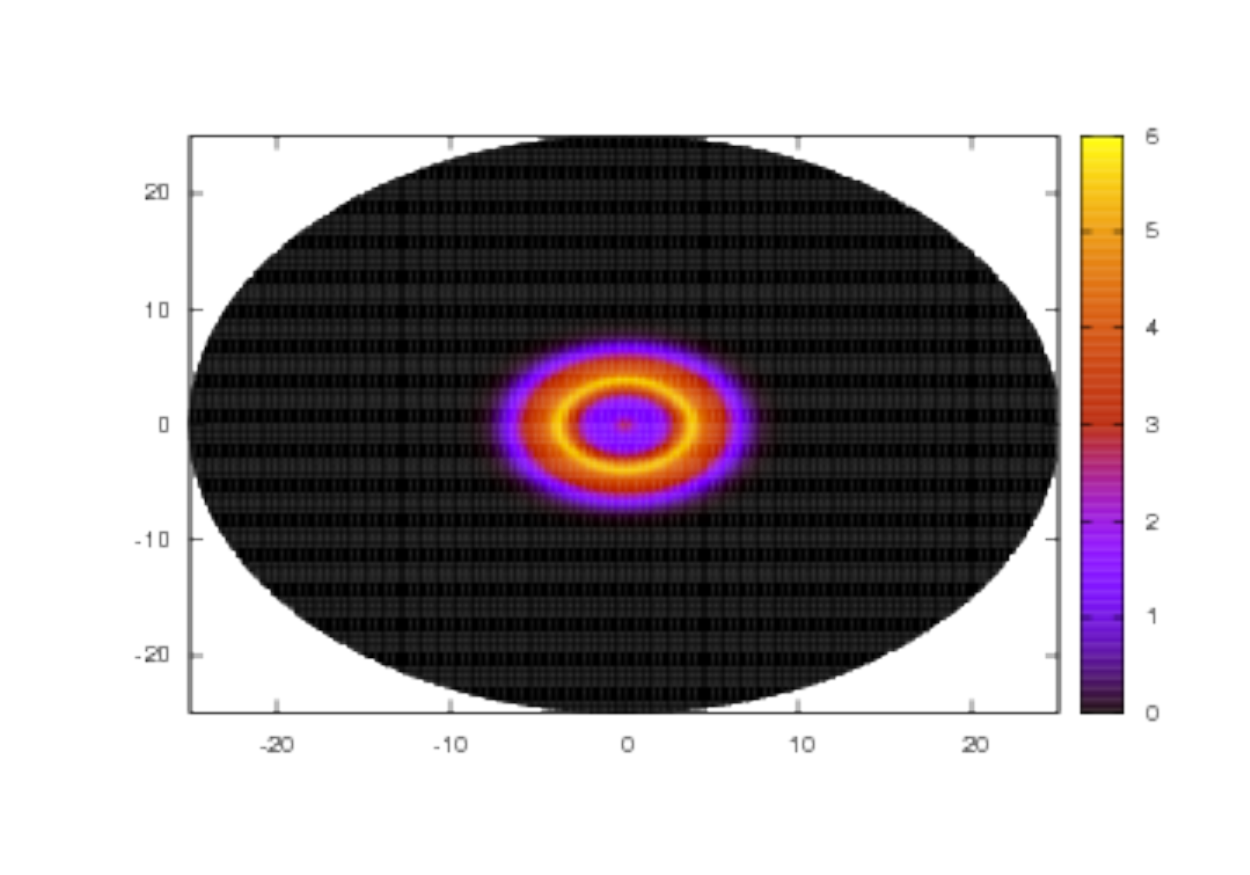} &
\includegraphics[width=5.cm,height=5.cm]{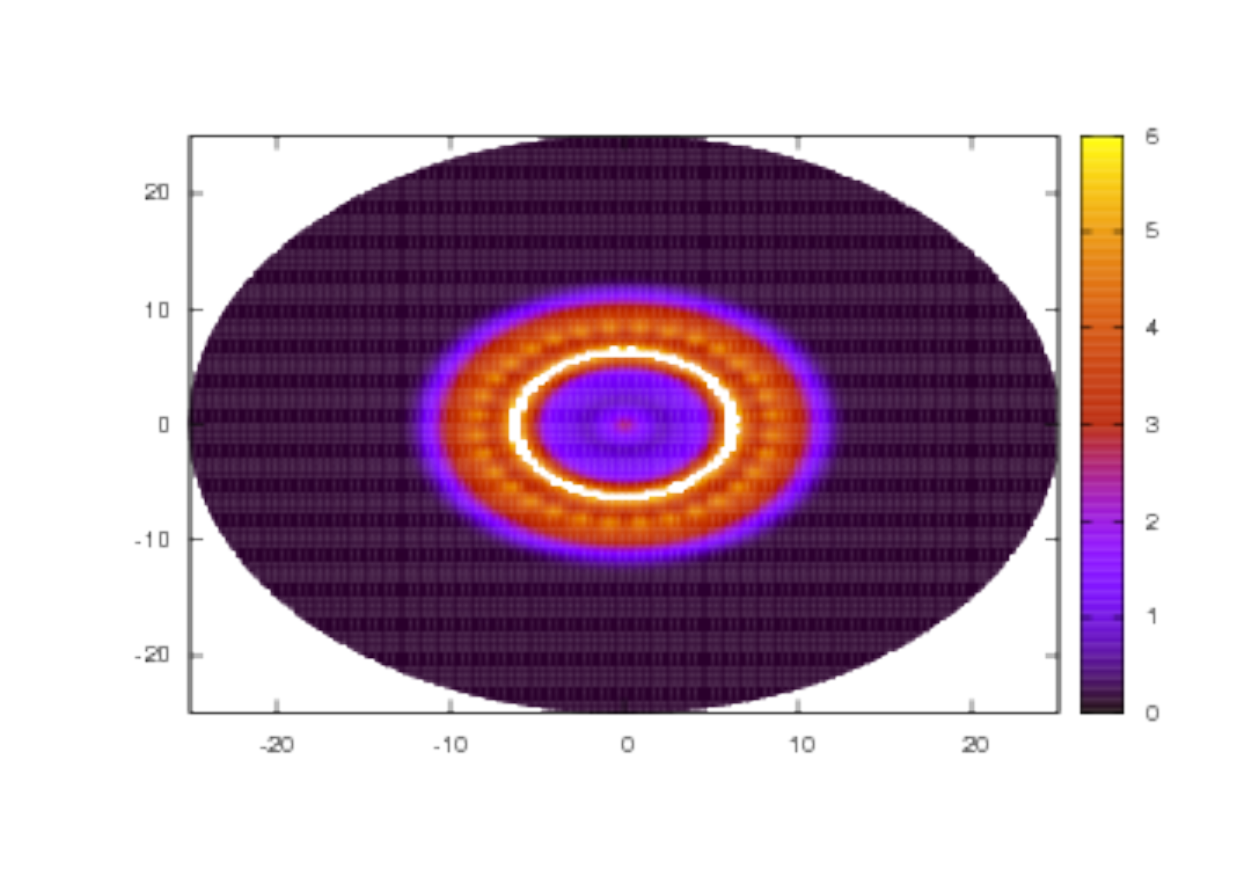} 
\\
\includegraphics[width=5.cm,height=5.cm]{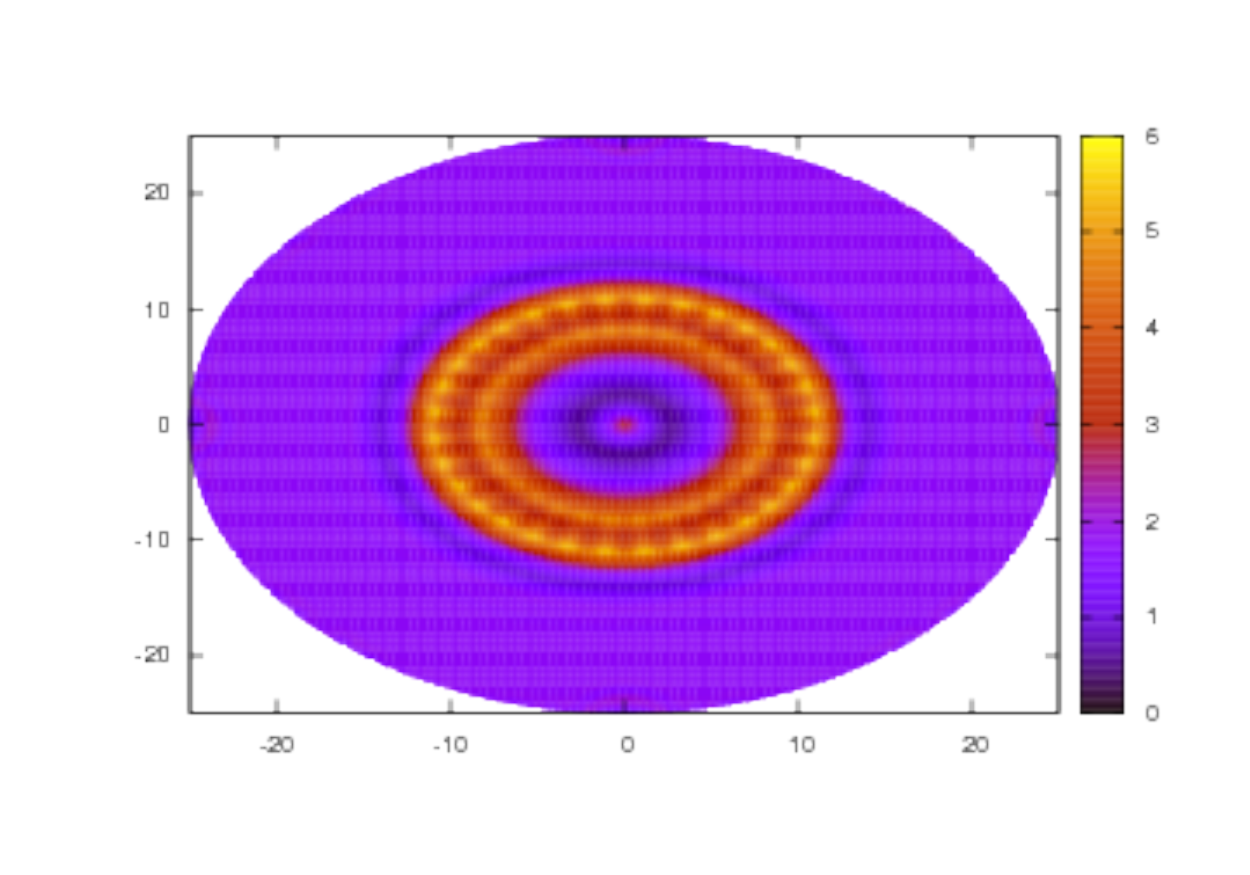} &
\includegraphics[width=5.cm,height=5.cm]{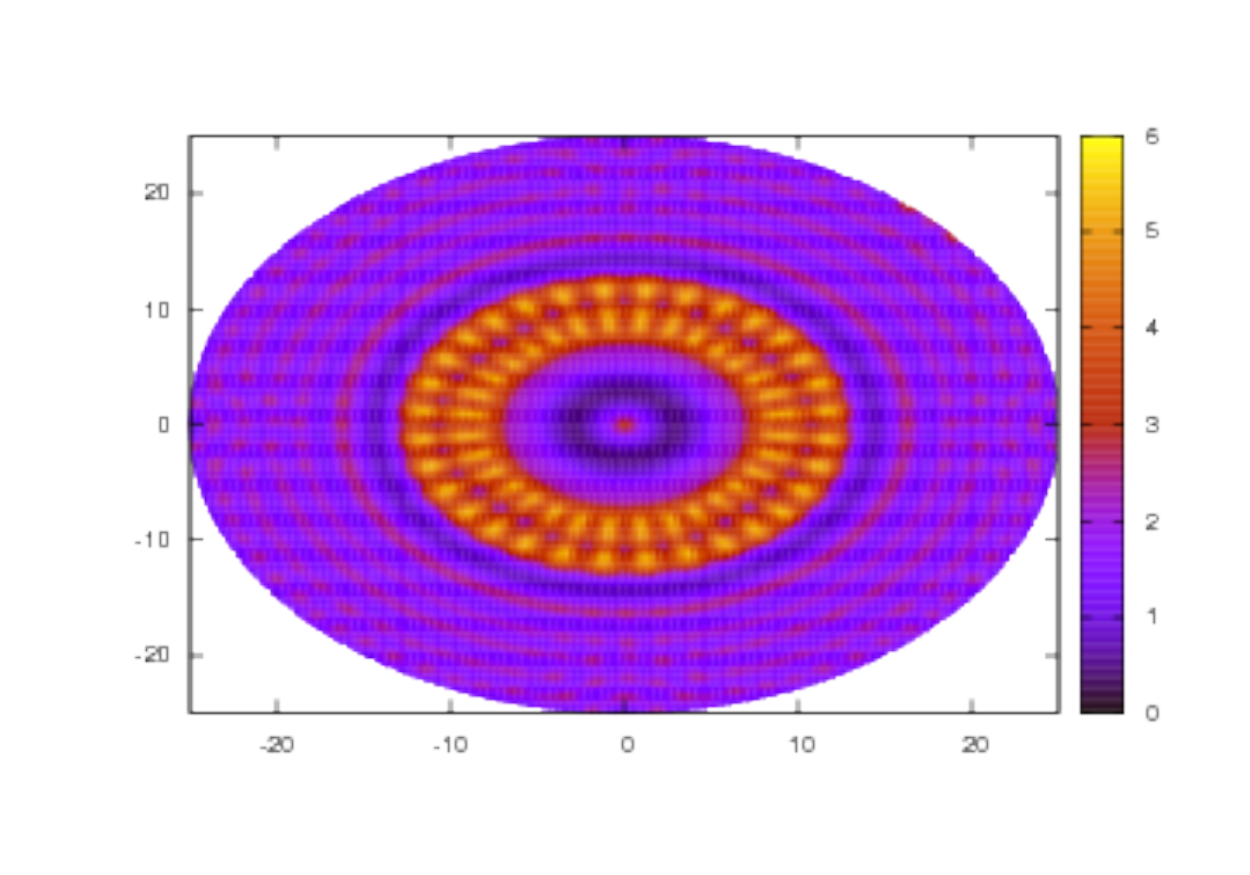} &
\includegraphics[width=5.cm,height=5.cm]{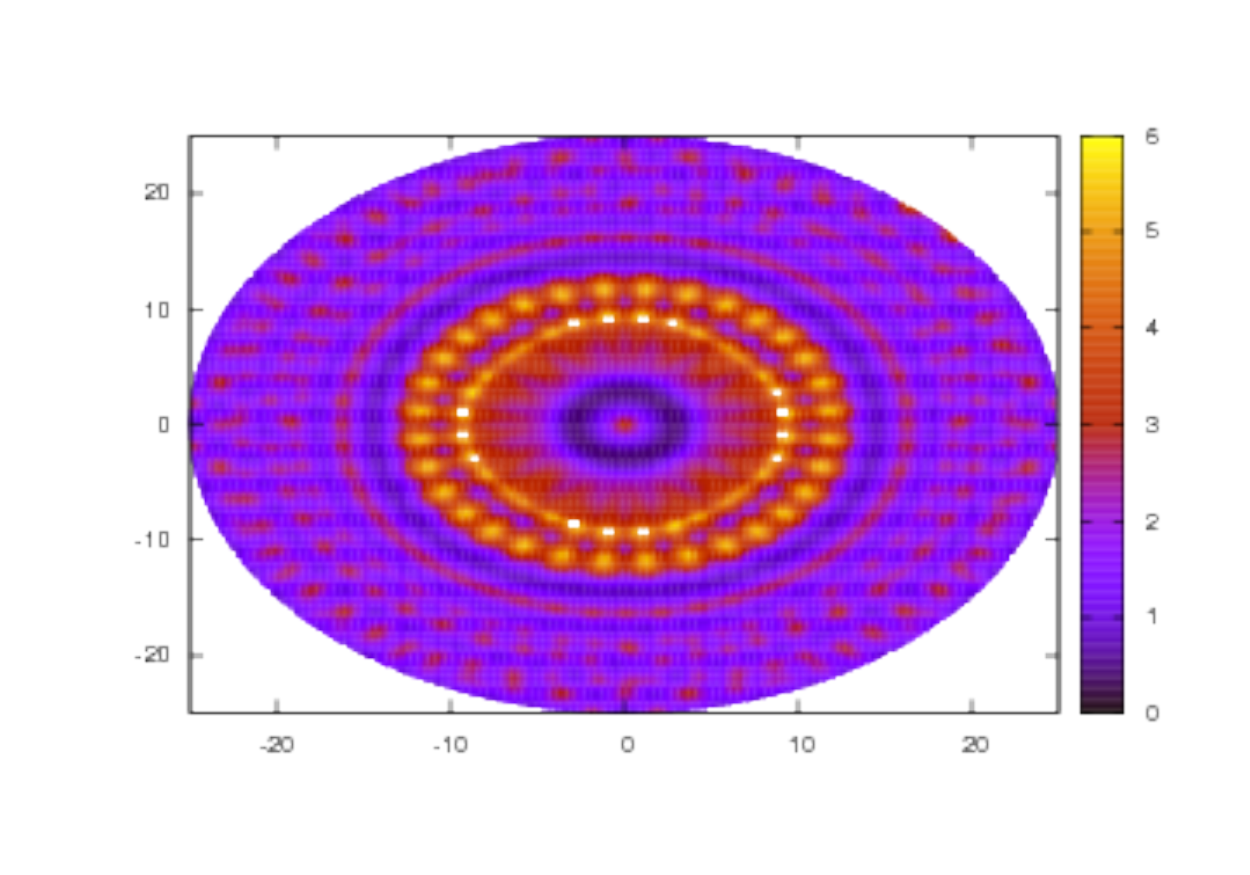} 
\end{tabular}
\caption{Test 5 : {\it Time evolution of the cell density (a) t = 0 $\overline{t}$ (b) t = 15 $\overline{t}$ (c) t = 25 $\overline{t}$ (d) t = 35 $\overline{t}$ (e) t = 45 $\overline{t}$ and (f) t = 55 $\overline{t}$ in sec.}}
\label{Fig4-1}
\end{figure}

\section{Conclusion}
\label{sec:conc}
\setcounter{equation}{0}

In this paper we present a new algorithm based on a Cartesian mesh for the numerical approximation of kinetic models on an arbitrary geometry boundary modelling chemosensitive movements. We present first a kinetic model for chemotactic bacteria  interacting with two chemical substances, {\it i.e.} nutrient and chemottractant. Then we give the numerical discretization for this kinetic model and the numerical method for the boundary conditions based on a Cartesian mesh. A large various numerical tests in $2D_x\times1D_v $ are shown and compared with biological experiences. We conclude that on the one hand this kinetic model  represents well the chemotactic bacteria behaviors and on the other hand our numerical method is accurate and efficient for numerical simulation.


  
\bibliographystyle{plain}

\end{document}